\documentclass[12pt]{amsart}      
\usepackage{amssymb}
\usepackage{eucal}
\usepackage{amsmath}
\usepackage{amscd}
\usepackage[dvips]{color}
\usepackage{multicol}
\usepackage[all]{xy}           
\usepackage{graphicx}
\usepackage{color}
\usepackage{colordvi}
\usepackage{xspace}
\usepackage{axodraw}

\begin{document}  

\newcommand{\nc}{\newcommand}
\newcommand{\delete}[1]{}
\nc{\dfootnote}[1]{{}}          
\nc{\ffootnote}[1]{\dfootnote{#1}}
\nc{\mfootnote}[1]{\footnote{#1}} 
\nc{\ofootnote}[1]{\footnote{\tiny Older version: #1}} 

\nc{\mlabel}[1]{\label{#1}}  
\nc{\mcite}[1]{\cite{#1}}  
\nc{\mref}[1]{\ref{#1}}  

\delete{
\nc{\mlabel}[1]{\label{#1}  
{\hfill \hspace{1cm}{\bf{{\ }\hfill(#1)}}}}
\nc{\mcite}[1]{\cite{#1}{{\bf{{\ }(#1)}}}}  
\nc{\mref}[1]{\ref{#1}{{\bf{{\ }(#1)}}}}  
}

\nc{\mbibitem}[1]{\bibitem{#1}} 
\nc{\mkeep}[1]{\marginpar{{\bf #1}}} 

\nc{\redtext}[1]{\textcolor{red}{#1}}
\nc{\bluetext}[1]{\textcolor{blue}{#1}}
\nc{\greentext}[1]{\textcolor{green}{#1}}

\nc{\li}[1]{\redtext{#1}}
\nc{\bill}[1]{\bluetext{#1}}

\newtheorem{theorem}{Theorem}[section]
\newtheorem{prop}[theorem]{Proposition}
\newtheorem{defn}[theorem]{Definition}
\newtheorem{lemma}[theorem]{Lemma}
\newtheorem{coro}[theorem]{Corollary}
\newtheorem{prop-def}{Proposition-Definition}[section]
\newtheorem{claim}{Claim}[section]
\newtheorem{remark}[theorem]{Remark}
\newtheorem{propprop}{Proposed Proposition}[section]
\newtheorem{conjecture}{Conjecture}
\newtheorem{exam}[theorem]{Example}
\newtheorem{assumption}{Assumption}
\newtheorem{condition}[theorem]{Assumption}
\newtheorem{question}[theorem]{Question}

\renewcommand{\labelenumi}{{\rm(\alph{enumi})}}
\renewcommand{\theenumi}{\alph{enumi}}

\nc{\tred}[1]{\textcolor{red}{#1}}
\nc{\tblue}[1]{\textcolor{blue}{#1}}
\nc{\tgreen}[1]{\textcolor{green}{#1}}
\nc{\tpurple}[1]{\textcolor{purple}{#1}}
\nc{\btred}[1]{\textcolor{red}{\bf #1}}
\nc{\btblue}[1]{\textcolor{blue}{\bf #1}}
\nc{\btgreen}[1]{\textcolor{green}{\bf #1}}
\nc{\btpurple}[1]{\textcolor{purple}{\bf #1}}

\nc{\adec}{\check{;}}
\nc{\dftimes}{\widetilde{\otimes}} \nc{\dfl}{\succ}
\nc{\dfr}{\prec} \nc{\dfc}{\circ} \nc{\dfb}{\bullet}
\nc{\dft}{\star} \nc{\dfcf}{{\mathbf k}} \nc{\spr}{\cdot}
\nc{\disp}[1]{\displaystyle{#1}}
\nc{\bin}[2]{ (_{\stackrel{\scs{#1}}{\scs{#2}}})}  
\nc{\binc}[2]{ \left (\!\! \begin{array}{c} \scs{#1}\\
    \scs{#2} \end{array}\!\! \right )}  
\nc{\bincc}[2]{  \left ( {\scs{#1} \atop
    \vspace{-.5cm}\scs{#2}} \right )}  
\nc{\sarray}[2]{\begin{array}{c}#1 \vspace{.1cm}\\ \hline
    \vspace{-.35cm} \\ #2 \end{array}}
\nc{\bs}{\bar{S}} \nc{\dcup}{\stackrel{\bullet}{\cup}}
\nc{\dbigcup}{\stackrel{\bullet}{\bigcup}} \nc{\etree}{\big |}
\nc{\la}{\longrightarrow} \nc{\fe}{\'{e}} \nc{\rar}{\rightarrow}
\nc{\dar}{\downarrow} \nc{\dap}[1]{\downarrow
\rlap{$\scriptstyle{#1}$}} \nc{\uap}[1]{\uparrow
\rlap{$\scriptstyle{#1}$}} \nc{\defeq}{\stackrel{\rm def}{=}}
\nc{\diffa}[1]{\{#1\}}
\nc{\diffs}[1]{\Delta{#1}}
\nc{\dis}[1]{\displaystyle{#1}} \nc{\dotcup}{\,
\displaystyle{\bigcup^\bullet}\ } \nc{\sdotcup}{\tiny{
\displaystyle{\bigcup^\bullet}\ }} \nc{\hcm}{\ \hat{,}\ }
\nc{\hcirc}{\hat{\circ}} \nc{\hts}{\hat{\shpr}}
\nc{\lts}{\stackrel{\leftarrow}{\shpr}}
\nc{\rts}{\stackrel{\rightarrow}{\shpr}} \nc{\lleft}{[}
\nc{\lright}{]} \nc{\uni}[1]{\tilde{#1}} \nc{\wor}[1]{\check{#1}}
\nc{\free}[1]{\bar{#1}} \nc{\den}[1]{\check{#1}} \nc{\lrpa}{\wr}
\nc{\curlyl}{\left \{ \begin{array}{c} {} \\ {} \end{array}
    \right .  \!\!\!\!\!\!\!}
\nc{\curlyr}{ \!\!\!\!\!\!\!
    \left . \begin{array}{c} {} \\ {} \end{array}
    \right \} }
\nc{\leaf}{\ell}       
\nc{\longmid}{\left | \begin{array}{c} {} \\ {} \end{array}
    \right . \!\!\!\!\!\!\!}
\nc{\ot}{\otimes} \nc{\sot}{{\scriptstyle{\ot}}}
\nc{\otm}{\overline{\ot}}
\nc{\ora}[1]{\stackrel{#1}{\rar}}
\nc{\ola}[1]{\stackrel{#1}{\la}}
\nc{\scs}[1]{\scriptstyle{#1}} \nc{\mrm}[1]{{\rm #1}}
\nc{\margin}[1]{\marginpar{\rm #1}}   
\nc{\dirlim}{\displaystyle{\lim_{\longrightarrow}}\,}
\nc{\invlim}{\displaystyle{\lim_{\longleftarrow}}\,}
\nc{\mvp}{\vspace{0.5cm}} \nc{\svp}{\vspace{2cm}}
\nc{\vp}{\vspace{8cm}} \nc{\proofbegin}{\noindent{\bf Proof: }}
\nc{\proofend}{$\blacksquare$ \vspace{0.5cm}}
\nc{\sha}{{\mbox{\cyr X}}}  
\nc{\ncsha}{{\mbox{\cyr X}^{\mathrm NC}}} \nc{\ncshao}{{\mbox{\cyr
X}^{\mathrm NC,\,0}}}
\nc{\shpr}{\diamond}    
\nc{\shprm}{\overline{\diamond}}    
\nc{\shpro}{\diamond^0}    
\nc{\shprr}{\diamond^r}     
\nc{\shpra}{\overline{\diamond}^r}
\nc{\shpru}{\check{\diamond}} \nc{\catpr}{\diamond_l}
\nc{\rcatpr}{\diamond_r} \nc{\lapr}{\diamond_a}
\nc{\sqcupm}{\ot}
\nc{\lepr}{\diamond_e} \nc{\vep}{\varepsilon} \nc{\labs}{\mid\!}
\nc{\rabs}{\!\mid} \nc{\hsha}{\widehat{\sha}}
\nc{\lsha}{\stackrel{\leftarrow}{\sha}}
\nc{\rsha}{\stackrel{\rightarrow}{\sha}} \nc{\lc}{\lfloor}
\nc{\rc}{\rfloor} \nc{\sqmon}[1]{\langle #1\rangle}
\nc{\forest}{\calf} \nc{\ass}[1]{\alpha({#1})}
\nc{\altx}{\Lambda_X} \nc{\vecT}{\vec{T}} \nc{\onetree}{\bullet}
\nc{\Ao}{\check{A}}
\nc{\seta}{\underline{\Ao}}
\nc{\deltaa}{\overline{\delta}}
\nc{\trho}{\tilde{\rho}}
\nc{\tpow}[2]{{#2}^{\ot #1}}

\nc{\mmbox}[1]{\mbox{\ #1\ }} \nc{\ann}{\mrm{ann}}
\nc{\Aut}{\mrm{Aut}}
\nc{\bread}{\mrm{b}}
\nc{\can}{\mrm{can}} \nc{\colim}{\mrm{colim}}
\nc{\Cont}{\mrm{Cont}} \nc{\rchar}{\mrm{char}}
\nc{\cok}{\mrm{coker}} \nc{\dtf}{{R-{\rm tf}}} \nc{\dtor}{{R-{\rm
tor}}}
\renewcommand{\det}{\mrm{det}}
\nc{\depth}{{\mrm d}}
\nc{\Div}{{\mrm Div}} \nc{\End}{\mrm{End}} \nc{\Ext}{\mrm{Ext}}
\nc{\Fil}{\mrm{Fil}} \nc{\Frob}{\mrm{Frob}} \nc{\Gal}{\mrm{Gal}}
\nc{\GL}{\mrm{GL}} \nc{\Hom}{\mrm{Hom}} \nc{\hsr}{\mrm{H}}
\nc{\hpol}{\mrm{HP}} \nc{\id}{\mrm{id}} \nc{\im}{\mrm{im}}
\nc{\incl}{\mrm{incl}} \nc{\length}{\mrm{length}}
\nc{\LR}{\mrm{LR}} \nc{\mchar}{\rm char} \nc{\NC}{\mrm{NC}}
\nc{\mpart}{\mrm{part}} \nc{\pl}{\mrm{PL}}
\nc{\ql}{{\QQ_\ell}} \nc{\qp}{{\QQ_p}}
\nc{\rank}{\mrm{rank}} \nc{\rba}{\rm{RBA }} \nc{\rbas}{\rm{RBAs }}
\nc{\rbpl}{\mrm{RBPL}}
\nc{\rbw}{\rm{RBW }} \nc{\rbws}{\rm{RBWs }} \nc{\rcot}{\mrm{cot}}
\nc{\rest}{\rm{controlled}\xspace}
\nc{\rdef}{\mrm{def}} \nc{\rdiv}{{\rm div}} \nc{\rtf}{{\rm tf}}
\nc{\rtor}{{\rm tor}} \nc{\res}{\mrm{res}} \nc{\SL}{\mrm{SL}}
\nc{\Spec}{\mrm{Spec}} \nc{\tor}{\mrm{tor}} \nc{\Tr}{\mrm{Tr}}
\nc{\mtr}{\mrm{sk}}

\nc{\ab}{\mathbf{Ab}} \nc{\Alg}{\mathbf{Alg}}
\nc{\Algo}{\mathbf{Alg}^0} \nc{\Bax}{\mathbf{Bax}}
\nc{\Baxo}{\mathbf{Bax}^0} \nc{\Dif}{\mathbf{Dif}}
\nc{\CDif}{\mathbf{CDif}}
\nc{\CRB}{\mathbf{CRB}}
\nc{\CDRB}{\mathbf{CDRB}}
\nc{\RB}{\mathbf{RB}}
\nc{\DRB}{\mathbf{DRB}}
\nc{\RBo}{\mathbf{RB}^0} \nc{\BRB}{\mathbf{RB}}
\nc{\Dend}{\mathbf{DD}}
\nc{\Set}{\mathbf{Set}}
\nc{\bfk}{{\bf k}} \nc{\bfone}{{\bf 1}}
\nc{\base}[1]{{a_{#1}}} \nc{\detail}{\marginpar{\bf More detail}
    \noindent{\bf Need more detail!}
    \svp}
\nc{\Diff}{\mathbf{Diff}} \nc{\gap}{\marginpar{\bf
Incomplete}\noindent{\bf Incomplete!!}
    \svp}
\nc{\FMod}{\mathbf{FMod}} \nc{\mset}{\mathbf{MSet}}
\nc{\rb}{\mathrm{RB}} \nc{\Int}{\mathbf{Int}}
\nc{\Mon}{\mathbf{Mon}}
\nc{\remarks}{\noindent{\bf Remarks: }} \nc{\Rep}{\mathbf{Rep}}
\nc{\Rings}{\mathbf{Rings}} \nc{\Sets}{\mathbf{Sets}}
\nc{\DT}{\mathbf{DT}}

\nc{\BA}{{\mathbb A}} \nc{\CC}{{\mathbb C}} \nc{\DD}{{\mathbb D}}
\nc{\EE}{{\mathbb E}} \nc{\FF}{{\mathbb F}} \nc{\GG}{{\mathbb G}}
\nc{\HH}{{\mathbb H}} \nc{\LL}{{\mathbb L}} \nc{\NN}{{\mathbb N}}
\nc{\QQ}{{\mathbb Q}} \nc{\RR}{{\mathbb R}} \nc{\TT}{{\mathbb T}}
\nc{\VV}{{\mathbb V}} \nc{\ZZ}{{\mathbb Z}}


\nc{\cala}{{\mathcal A}} \nc{\calc}{{\mathcal C}}
\nc{\cald}{{\mathcal D}} \nc{\cale}{{\mathcal E}}
\nc{\calf}{{\mathcal F}} \nc{\calfr}{{{\mathcal F}^{\,r}}}
\nc{\calfo}{{\mathcal F}^0} \nc{\calfro}{{\mathcal F}^{\,r,0}}
\nc{\oF}{\overline{F}}  \nc{\calg}{{\mathcal G}}
\nc{\calh}{{\mathcal H}} \nc{\cali}{{\mathcal I}}
\nc{\calj}{{\mathcal J}} \nc{\call}{{\mathcal L}}
\nc{\calm}{{\mathcal M}} \nc{\caln}{{\mathcal N}}
\nc{\calo}{{\mathcal O}} \nc{\calp}{{\mathcal P}}
\nc{\calr}{{\mathcal R}} \nc{\calt}{{\mathcal T}}
\nc{\caltr}{{\mathcal T}^{\,r}}
\nc{\calu}{{\mathcal U}} \nc{\calv}{{\mathcal V}}
\nc{\calw}{{\mathcal W}} \nc{\calx}{{\mathcal X}}
\nc{\CA}{\mathcal{A}}

\nc{\fraka}{{\mathfrak a}} \nc{\frakB}{{\mathfrak B}}
\nc{\frakb}{{\mathfrak b}} \nc{\frakd}{{\mathfrak d}}
\nc{\oD}{\overline{D}}
\nc{\frakF}{{\mathfrak F}} \nc{\frakg}{{\mathfrak g}}
\nc{\frakm}{{\mathfrak m}} \nc{\frakM}{{\mathfrak M}}
\nc{\frakMo}{{\mathfrak M}^0} \nc{\frakp}{{\mathfrak p}}
\nc{\frakS}{{\mathfrak S}} \nc{\frakSo}{{\mathfrak S}^0}
\nc{\fraks}{{\mathfrak s}} \nc{\os}{\overline{\fraks}}
\nc{\frakT}{{\mathfrak T}}
\nc{\oT}{\overline{T}}
\nc{\frakX}{{\mathfrak X}} \nc{\frakXo}{{\mathfrak X}^0}
\nc{\frakx}{{\mathbf x}}
\nc{\frakTx}{\frakT}      
\nc{\frakTa}{\frakT^a}        
\nc{\frakTxo}{\frakTx^0}   
\nc{\caltao}{\calt^{a,0}}   
\nc{\oV}{\overline{V}}
\nc{\ox}{\overline{\frakx}} \nc{\fraky}{{\mathfrak y}}
\nc{\frakz}{{\mathfrak z}} \nc{\oX}{\overline{X}}
\nc{\oZ}{\overline{Z}}

\font\cyr=wncyr10


\def\ta1{{\scalebox{0.25}{ 
\begin{picture}(12,12)(38,-38)
\SetWidth{0.5} \SetColor{Black} \Vertex(45,-33){5.66}
\end{picture}}}}

\def\tb2{{\scalebox{0.25}{ 
\begin{picture}(12,42)(38,-38)
\SetWidth{0.5} \SetColor{Black} \Vertex(45,-3){5.66}
\SetWidth{1.0} \Line(45,-3)(45,-33) \SetWidth{0.5}
\Vertex(45,-33){5.66}
\end{picture}}}}

\def\tc3{{\scalebox{0.25}{ 
\begin{picture}(12,72)(38,-38)
\SetWidth{0.5} \SetColor{Black} \Vertex(45,27){5.66}
\SetWidth{1.0} \Line(45,27)(45,-3) \SetWidth{0.5}
\Vertex(45,-33){5.66} \SetWidth{1.0} \Line(45,-3)(45,-33)
\SetWidth{0.5} \Vertex(45,-3){5.66}
\end{picture}}}}

\def\td31{{\scalebox{0.25}{ 
\begin{picture}(42,42)(23,-38)
\SetWidth{0.5} \SetColor{Black} \Vertex(45,-3){5.66}
\Vertex(30,-33){5.66} \Vertex(60,-33){5.66} \SetWidth{1.0}
\Line(45,-3)(30,-33) \Line(60,-33)(45,-3)
\end{picture}}}}

\def\te4{{\scalebox{0.25}{ 
\begin{picture}(12,102)(38,-8)
\SetWidth{0.5} \SetColor{Black} \Vertex(45,57){5.66}
\Vertex(45,-3){5.66} \Vertex(45,27){5.66} \Vertex(45,87){5.66}
\SetWidth{1.0} \Line(45,57)(45,27) \Line(45,-3)(45,27)
\Line(45,57)(45,87)
\end{picture}}}}

\def\tf41{{\scalebox{0.25}{ 
\begin{picture}(42,72)(38,-8)
\SetWidth{0.5} \SetColor{Black} \Vertex(45,27){5.66}
\Vertex(45,-3){5.66} \SetWidth{1.0} \Line(45,27)(45,-3)
\SetWidth{0.5} \Vertex(60,57){5.66} \SetWidth{1.0}
\Line(45,27)(60,57) \SetWidth{0.5} \Vertex(75,27){5.66}
\SetWidth{1.0} \Line(75,27)(60,57)
\end{picture}}}}

\def\tg42{{\scalebox{0.25}{ 
\begin{picture}(42,72)(8,-8)
\SetWidth{0.5} \SetColor{Black} \Vertex(45,27){5.66}
\Vertex(45,-3){5.66} \SetWidth{1.0} \Line(45,27)(45,-3)
\SetWidth{0.5} \Vertex(15,27){5.66} \Vertex(30,57){5.66}
\SetWidth{1.0} \Line(15,27)(30,57) \Line(45,27)(30,57)
\end{picture}}}}

\def\th43{{\scalebox{0.25}{ 
\begin{picture}(42,42)(8,-8)
\SetWidth{0.5} \SetColor{Black} \Vertex(45,-3){5.66}
\Vertex(15,-3){5.66} \Vertex(30,27){5.66} \SetWidth{1.0}
\Line(15,-3)(30,27) \Line(45,-3)(30,27) \Line(30,27)(30,-3)
\SetWidth{0.5} \Vertex(30,-3){5.66}
\end{picture}}}}

\def\thII43{{\scalebox{0.25}{ 
\begin{picture}(72,57) (68,-128)
    \SetWidth{0.5}
    \SetColor{Black}
    \Vertex(105,-78){5.66}
    \SetWidth{1.5}
    \Line(105,-78)(75,-123)
    \Line(105,-78)(105,-123)
    \Line(105,-78)(135,-123)
    \SetWidth{0.5}
    \Vertex(75,-123){5.66}
    \Vertex(105,-123){5.66}
    \Vertex(135,-123){5.66}
  \end{picture}
  }}}

\def\thj44{{\scalebox{0.25}{ 
\begin{picture}(42,72)(8,-8)
\SetWidth{0.5} \SetColor{Black} \Vertex(30,57){5.66}
\SetWidth{1.0} \Line(30,57)(30,27) \SetWidth{0.5}
\Vertex(30,27){5.66} \SetWidth{1.0} \Line(45,-3)(30,27)
\SetWidth{0.5} \Vertex(45,-3){5.66} \Vertex(15,-3){5.66}
\SetWidth{1.0} \Line(15,-3)(30,27)
\end{picture}}}}

\def\ti5{{\scalebox{0.25}{ 
\begin{picture}(12,132)(23,-8)
\SetWidth{0.5} \SetColor{Black} \Vertex(30,117){5.66}
\SetWidth{1.0} \Line(30,117)(30,87) \SetWidth{0.5}
\Vertex(30,87){5.66} \Vertex(30,57){5.66} \Vertex(30,27){5.66}
\Vertex(30,-3){5.66} \SetWidth{1.0} \Line(30,-3)(30,27)
\Line(30,27)(30,57) \Line(30,87)(30,57)
\end{picture}}}}

\def\tj51{{\scalebox{0.25}{ 
\begin{picture}(42,102)(53,-38)
\SetWidth{0.5} \SetColor{Black} \Vertex(61,27){4.24}
\SetWidth{1.0} \Line(75,57)(90,27) \Line(60,27)(75,57)
\SetWidth{0.5} \Vertex(90,-3){5.66} \Vertex(60,27){5.66}
\Vertex(75,57){5.66} \Vertex(90,-33){5.66} \SetWidth{1.0}
\Line(90,-33)(90,-3) \Line(90,-3)(90,27) \SetWidth{0.5}
\Vertex(90,27){5.66}
\end{picture}}}}

\def\tk52{{\scalebox{0.25}{ 
\begin{picture}(42,102)(23,-8)
\SetWidth{0.5} \SetColor{Black} \Vertex(60,57){5.66}
\Vertex(45,87){5.66} \SetWidth{1.0} \Line(45,87)(60,57)
\SetWidth{0.5} \Vertex(30,57){5.66} \SetWidth{1.0}
\Line(30,57)(45,87) \SetWidth{0.5} \Vertex(30,-3){5.66}
\SetWidth{1.0} \Line(30,-3)(30,27) \SetWidth{0.5}
\Vertex(30,27){5.66} \SetWidth{1.0} \Line(30,57)(30,27)
\end{picture}}}}

\def\tl53{{\scalebox{0.25}{ 
\begin{picture}(42,102)(8,-8)
\SetWidth{0.5} \SetColor{Black} \Vertex(30,57){5.66}
\Vertex(30,27){5.66} \SetWidth{1.0} \Line(30,57)(30,27)
\SetWidth{0.5} \Vertex(30,87){5.66} \SetWidth{1.0}
\Line(30,27)(45,-3) \SetWidth{0.5} \Vertex(15,-3){5.66}
\SetWidth{1.0} \Line(15,-3)(30,27) \Line(30,57)(30,87)
\SetWidth{0.5} \Vertex(45,-3){5.66}
\end{picture}}}}

\def\tm54{{\scalebox{0.25}{ 
\begin{picture}(42,72)(8,-38)
\SetWidth{0.5} \SetColor{Black} \Vertex(30,-3){5.66}
\SetWidth{1.0} \Line(30,27)(30,-3) \Line(30,-3)(45,-33)
\SetWidth{0.5} \Vertex(15,-33){5.66} \SetWidth{1.0}
\Line(15,-33)(30,-3) \SetWidth{0.5} \Vertex(45,-33){5.66}
\SetWidth{1.0} \Line(30,-33)(30,-3) \SetWidth{0.5}
\Vertex(30,-33){5.66} \Vertex(30,27){5.66}
\end{picture}}}}

\def\tn55{{\scalebox{0.25}{ 
\begin{picture}(42,72)(8,-38)
\SetWidth{0.5} \SetColor{Black} \Vertex(15,-33){5.66}
\Vertex(45,-33){5.66} \Vertex(30,27){5.66} \SetWidth{1.0}
\Line(45,-33)(45,-3) \SetWidth{0.5} \Vertex(45,-3){5.66}
\Vertex(15,-3){5.66} \SetWidth{1.0} \Line(30,27)(45,-3)
\Line(15,-3)(30,27) \Line(15,-3)(15,-33)
\end{picture}}}}


\def\ydec31{\!\!\includegraphics[scale=0.5]{ydec31.eps}}

\def\kyldec31{\!\!\includegraphics[scale=0.5]{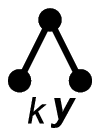}}

\def\yldec31{\!\!\includegraphics[scale=0.5]{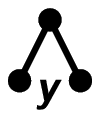}}

\def\xldec41r{\!\!\includegraphics[scale=0.5]{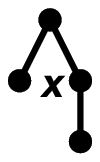}}

\def\xyldec43{\!\!\includegraphics[scale=0.5]{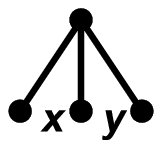}}

\def\xtd31{\!\!\includegraphics[scale=1]{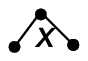}}

\def\xthj44{\!\!\includegraphics[scale=1]{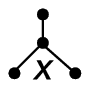}}


\title{On Differential Rota-Baxter algebras}
\author{Li Guo}
\address{Department of Mathematics and Computer Science,
         Rutgers University,
         Newark, NJ 07102}
\email{liguo@newark.rutgers.edu}
\author{William Keigher}
\address{Department of Mathematics and Computer Science,
         Rutgers University,
         Newark, NJ 07102}
\email{keigher@andromeda.rutgers.edu}



\begin{abstract}
A Rota-Baxter operator of weight $\lambda$ is an abstraction of both the integral operator (when $\lambda=0$) and the summation operator (when $\lambda=1$). We similarly define a differential operator of weight $\lambda$ that includes both the differential operator (when $\lambda=0$) and the difference operator (when $\lambda=1$).
We further consider an algebraic structure with both a differential operator of weight $\lambda$
and a Rota-Baxter operator of weight $\lambda$ that are related in the same way that the differential operator and the integral operator are related by the First Fundamental Theorem of Calculus. We construct free objects in the corresponding categories. In the commutative case, the free objects are given in terms of generalized shuffles, called mixable shuffles. In the noncommutative case, the free objects are given in terms of angularly decorated rooted forests.
As a byproduct, we obtain structures of a differential algebra
on decorated and undecorated planar rooted forests.
\end{abstract}


\maketitle

\noindent
{\bf Keywords: } Differential algebra of weight $\lambda$, Rota-Baxter algebra, differential-Rota-Baxter algebra, construction of free algebras, shuffle product, planar trees.


\setcounter{section}{0}

\section{Introduction}
\subsection{Motivation}
The First Fundamental Theorem of Calculus states that
(under suitable conditions)
\begin{equation}
\frac{d}{dx} \Big(\int_a^x f(t)dt\Big) = f(x).
\mlabel{eq:cal}
\end{equation}
Thus the integral operator $P(f)(x)=\int_a^x f(t)dt$ is
the right inverse of the differential operator $d(f)(x)=\frac{d f}{dx}(x)$, so that $(d \circ P)(f)=f$.
A similar relation holds for the difference operator and summation operator (see Example~\mref{ex:diff}.(\mref{it:lsum})). The abstraction of the differential operator and difference operator led to the development of differential algebra and difference algebra~\mcite{RCo,Kol}. Likewise, the integral operator $P$ and summation operator have been abstracted to give the notion of Rota-Baxter operators (previously called Baxter operators) and Rota-Baxter algebras~\mcite{Ba,Rot1,Rot3}.
In the last few years, major progresses have been made in both differential algebra and Rota-Baxter algebra, with applications in broad areas in mathematics and physics~\mcite{Ag3,A-M,DART,C-K1,E-G4,E-G0,EGK3,EGM,G-K1,S-P1,S-P}. For instance, both operators played important roles in the recent developments in renormalization of quantum field
theory~\mcite{C-K1,C-M,E-Gs}.

This paper studies the algebraic structure reflecting the relation between the differential operator and the integral operator as in the First Fundamental Theorem of Calculus. By analogy to a Rota-Baxter operator that unifies the notions of an integral operator and a summation operator, we first unify the concepts of the differential operator and the difference operator by the concept of a $\lambda$-differential operator, where $\lambda$ is a fixed element in the ground ring, that gives the differential (resp. difference) operator when $\lambda$ is 0 (resp. 1). We then introduce the concept of a differential Rota-Baxter algebra of weight $\lambda$ consisting of an algebra with both a $\lambda$-differential operator and a $\lambda$-Rota-Baxter operator with a compatibility condition between these two operators.

\subsection{Definitions and preliminary examples}
\begin{defn} {\rm
Let $\bfk$ be a unitary commutative ring. Let $\lambda\in \bfk$
be fixed.
\begin{enumerate}
\item
A {\bf differential $\bfk$-algebra of weight $\lambda$} (also called a {\bf $\lambda$-differential $\bfk$-algebra}) is an associative $\bfk$-algebra $R$ together with a linear operator $d:R\to R$ such that
\begin{equation}
d(xy)=d(x)y+x d(y)+ \lambda d(x)d(y), \forall\, x,y\in R, \mlabel{eq:diff}
\end{equation}
and
\begin{equation}
d(1)=0.
\mlabel{eq:diffc}
\end{equation}
Such an operator is called a {\bf differential operator of weight $\lambda$} or a {\bf derivation of weight $\lambda$}. It is also called a {\bf $\lambda$-differential operator} or a {\bf $\lambda$-derivation}.
The category of differential algebras (resp. commutative differential algebras) of weight $\lambda$ is denoted by $\Dif_\lambda$ (resp. $\CDif_\lambda$).
\item
A {\bf Rota-Baxter $\bfk$-algebra of weight $\lambda$} is an associative $\bfk$-algebra $R$ together with a linear operator $P:R\to R$ such that
\begin{equation}
P(x)P(y)=P(xP(y))+P(P(x)y)+ \lambda P(xy), \forall\, x,y\in R.
\mlabel{eq:rba}
\end{equation}
Such an operator is called a {\bf Rota-Baxter operator of weight $\lambda$} or a {\bf $\lambda$-Rota-Baxter operator}.
The category of Rota-Baxter algebras (resp. commutative Rota-Baxter algebras) of weight $\lambda$ is denoted by $\RB_\lambda$ (resp. $\CRB_\lambda$).
\item
A {\bf differential Rota-Baxter $\bfk$-algebra of weight $\lambda$} (also called a {\bf $\lambda$-differential Rota-Baxter $\bfk$-algebra})  is an associative $\bfk$-algebra $R$ together with a differential operator $d$ of weight $\lambda$ and a Rota-Baxter operator $P$ of weight $\lambda$ such that
\begin{equation}
d\circ P= \id_R.
\mlabel{eq:Baxdiff}
\end{equation}
The category of differential Rota-Baxter algebras (resp. commutative differential Rota-Baxter algebras) of weight $\lambda$ is denoted by $\DRB_\lambda$ (resp. $\CDRB_\lambda$).
\end{enumerate}
}
\mlabel{de:main}
\end{defn}
We also use $\Alg=\Alg_\bfk$ to denote the category of $\bfk$-algebras. When there is no danger of confusion, we will suppress $\lambda$ and $\bfk$ from the notations.
We will also denote $\NN$ for the set of non-negative integers and $\NN_+$ for the set of positive integers.

Note that we require that a differential operator $d$ satisfies $d(1)=0$. A linear operator $d$ satisfying Eq.~(\mref{eq:diff}) is called a {\bf weak differential operator of weight $\lambda$}. A weak differential operator of weight $\lambda$ with $d(1)\neq 0$ is called a {\bf degenerated differential operator} of weight $\lambda$ for the reason given in Remark~\mref{re:deg}, and will be discussed in Section~\mref{ss:diff}.

We give some simple examples of differential, Rota-Baxter and differential Rota-Baxter algebras. Further examples will be given in later sections.

\begin{exam} {\rm
\begin{enumerate}
\item
A $0$-derivation and a $0$-differential algebra is a derivation and differential algebra in the usual sense~\mcite{Kol}.
\item
Let $\lambda \in \RR$,
$\lambda \neq 0$.  Let $R=\Cont(\RR)$ denote the $\RR$-algebra of continuous functions $f: \RR\to \RR$, and consider the usual
"difference quotient" operator $d_\lambda$ on $R$ defined by
\begin{equation}
(d_\lambda(f))(x) = (f(x+\lambda) - f(x))/\lambda.
\mlabel{eq:ldiff}
\end{equation}
Then it is immediate that $d_\lambda$ is a $\lambda$-derivation on $R$.
When $\lambda=1$, we obtain the usual difference operator on functions. Further, the usual derivation is
$\disp{d_0:= \lim_{\lambda \to 0} d_\lambda.}$
\item
A difference algebra~\mcite{RCo} is defined to be a commutative algebra $R$ together with an injective algebra endomorphism $\phi$ on $R$. It is simple to check that $\phi-\id$ is a differential operator of weight $1$.
\item
By the First Fundamental Theorem of Calculus in Eq.~(\mref{eq:cal}),
$(\Cont(\RR),d/dx, \int_0^x)$ is a differential Rota-Baxter algebra of weight 0.
\item
Let $0<\lambda \in \RR$. Let $R$ be an $\RR$-subalgebra of $\Cont(\RR)$ that is closed under the operators
$$P_0(f)(x)=-\int_x^\infty f(t)dt,\quad  P_\lambda(f)(x)=-\lambda\sum_{n\geq 0} f(x+n\lambda).$$
For example, $R$ can be taken to be the $\RR$-subalgebra generated by $e^{-x}$: $R=\sum_{k\geq 1} \RR e^{-kx}$.
Then $P_\lambda$ is a Rota-Baxter operator of weight $\lambda$ and, for the $d_\lambda$ in Eq.~(\mref{eq:ldiff}),
$$ d_\lambda\circ P_\lambda=\id_R, \forall\ 0\neq \lambda\in \RR, $$
reducing to the fundamental theorem
$ d_0\circ P_0=\id_R$ when $\lambda$ goes to $0$.
So $(R,d_\lambda,P_\lambda)$ is a differential Rota-Baxter algebra of weight $\lambda$.
\mlabel{it:lsum}
\end{enumerate}
}
\mlabel{ex:diff}
\end{exam}

\subsection{Main results and outline of the paper}
Our main purpose in this paper is to construct free objects in the various categories of $\lambda$-differential algebras and $\lambda$-differential Rota-Baxter algebras.

In Section~\mref{sec:diff}, we first prove basic properties of $\lambda$-differential algebras. We then construct the free objects in $\Dif_\lambda$ in Theorem~\mref{thm:freediff} and cofree objects in $\Dif_\lambda$ in Corollary~\mref{coro:cofree}. The construction of free objects in $\CDRB_\lambda$ is carried out in Section~\mref{sec:commfree} (Theorem~\mref{thm:commdiffrb}) and the construction of free objects in $\DRB_\lambda$ is carried out in Section~\mref{sec:ncfree} (Theorem~\mref{thm:ncdiffrb}). Both constructions rely on the explicit construction of free Rota-Baxter algebras,  in the commutative case in~\mcite{G-K1,G-K2} and in the noncommutative case in~\mcite{A-M,E-G4,E-G0}. Consequently, we obtain a structure of a differential algebra on the mixable shuffle and shuffle algebras, and on angularly decorated rooted trees. We further obtain the structure of a $\lambda$-differential algebra on planar rooted forests in Section~\mref{sec:tree} (Theorem~\mref{thm:treediff}).
It would be interesting to see how this is related to the work of Grossman and Larson~\mcite{G-L} on differential algebra structures on trees.
\medskip 

\noindent
{\bf Acknowledgements: } The first named author acknowledges support from NSF grant DMS-0505643.

\section{Differential algebras of weight $\lambda$}
\mlabel{sec:diff}
We first give some basic properties of $\lambda$-differential algebras, followed by a study of free and cofree $\lambda$-differential algebras.

\subsection{Basic properties and degenerated differential operators}
\mlabel{ss:diff}
Some basic properties of differential operators can be
easily generalized to $\lambda$-differential operators.
The following proposition generalizes the power rule in differential calculus and the well-known result
of Leibniz~\cite[p.60]{Kol}. It holds without the assumption that $d(1)=0$.

\begin{prop}
Let $(R,d)$ be a differential $\bfk$-algebra of weight $\lambda$.
\begin{enumerate}
\item
Let $x\in R$ and $n\in \NN_+$. Then
\[ d(x^n) =\sum_{i=1}^n \bincc{n}{i} \lambda^{i-1} x^{n-i} d(x)^i.\]
\mlabel{it:power}
\item
Let $x,y \in
R$, and let $n \in \NN$.
Then
\begin{equation}
d^{n}(xy) =
\sum_{k=0}^{n}\sum_{j=0}^{n-k}\binc{n}{k}\binc{n-k}{j}
\lambda^{k}d^{n-j}(x)d^{k+j}(y).
\label{eq:der2}
\end{equation}
\mlabel{it:prod}
\end{enumerate}
\mlabel{pp:basic}
\end{prop}
\begin{proof} (\mref{it:power}) The proof is similar to the inductive proof on $n$ for the usual power rule, using an index shift and Pascal's rule.
\smallskip

\noindent
(\mref{it:prod})
The proof is again similar to the case for differential operators. Proceeding by induction on $n$, the case $n = 0$ is trivial, so
assume that equation~(\ref{eq:der2}) holds for $n$, and consider
\begin{equation}
d^{n+1}(xy) = d^{n}(d(xy))
= d^{n}(d(x)y) + d^{n}(xd(y)) + \lambda d^{n}(d(x)d(y))).
\mlabel{eq:prodpf}
\end{equation}
Applying the induction hypothesis to the first term gives
\allowdisplaybreaks{
\begin{eqnarray*}
&&\sum_{k=0}^{n}\sum_{j=0}^{n-k}\binc{n}{k}\binc{n-k}{j}
\lambda^{k}d^{n+1-j}(x)d^{k+j}(y)\\
&=& \sum_{k=0}^{n}\sum_{j=1}^{n-k}\binc{n}{k}\binc{n-k}{j}
\lambda^{k}d^{n+1-j}(x)d^{k+j}(y)
+ \sum_{k=0}^{n}\binc{n}{k}
\lambda^{k}d^{n+1}(x)d^{k}(y).
\end{eqnarray*}
}
Doing the same to the second term in Eq.~(\mref{eq:prodpf}) followed by an index shift gives

\allowdisplaybreaks{
\begin{eqnarray*}
&&\sum_{k=0}^{n}\sum_{j=0}^{n-k}\binc{n}{k}\binc{n-k}{j}
\lambda^{k}d^{n-j}(x)d^{k+j+1}(y)\\
&=& \sum_{k=0}^{n}\sum_{j=1}^{n+1-k}\binc{n}{k}\binc{n-k}{j-1}
\lambda^{k}d^{n+1-j}(x)d^{k+j}(y)\\
&=& \sum_{k=0}^{n}\sum_{j=1}^{n-k}\binc{n}{k}\binc{n-k}{j-1}
\lambda^{k}d^{n+1-j}(x)d^{k+j}(y)
+ \sum_{k=0}^{n}\binc{n}{k}
\lambda^{k}d^{k}(x)d^{n+1}(y).
\end{eqnarray*}
}
Thus by Pascal's rule,
\allowdisplaybreaks{
\begin{eqnarray*}
d^{n}(d(x)y) + d^{n}(xd(y))&=& \sum_{k=0}^{n}\sum_{j=1}^{n-k}\binc{n}{k}\binc{n+1-k}{j}
\lambda^{k}d^{n+1-j}(x)d^{k+j}(y)\\
&& + \sum_{k=0}^{n}\binc{n}{k}
\lambda^{k}d^{n+1}(x)d^{k}(y)
+ \sum_{k=0}^{n}\binc{n}{k}
\lambda^{k}d^{k}(x)d^{n+1}(y)\\
&=& \sum_{k=0}^{n}\sum_{j=0}^{n+1-k}\binc{n}{k}\binc{n+1-k}{j}
\lambda^{k}d^{n+1-j}(x)d^{k+j}(y)\\
&=& \sum_{k=1}^{n}\sum_{j=0}^{n+1-k}\binc{n}{k}\binc{n+1-k}{j}
\lambda^{k}d^{n+1-j}(x)d^{k+j}(y)\\
&& + \sum_{j=0}^{n+1}\binc{n+1}{j}d^{n+1-j}(x)d^{j}(y).
\end{eqnarray*}
}
For the same reason, the third term in Eq.~(\mref{eq:prodpf}) gives
\allowdisplaybreaks{
\begin{eqnarray*}
&&\sum_{k=0}^{n}\sum_{j=0}^{n-k}\binc{n}{k}\binc{n-k}{j}
\lambda^{k+1}d^{n+1-j}(x)d^{k+j+1}(y)\\
&=& \sum_{k=1}^{n+1}\sum_{j=0}^{n+1-k}\binc{n}{k-1}\binc{n+1-k}{j}
\lambda^{k}d^{n+1-j}(x)d^{k+j}(y)\\
&=& \sum_{k=1}^{n}\sum_{j=0}^{n+1-k}\binc{n}{k-1}\binc{n+1-k}{j}
\lambda^{k}d^{n+1-j}(x)d^{k+j}(y)
+ \lambda^{n+1}d^{n+1}(x)d^{n+1}(y).
\end{eqnarray*}
}
Therefore another application of Pascal's rule gives
\allowdisplaybreaks{
\begin{eqnarray*}
&& d^{n}(d(x)y) + d^{n}(xd(y)) + \lambda d^{n}(d(x)d(y)))\\
&=& \sum_{k=1}^{n}\sum_{j=0}^{n+1-k}\binc{n+1}{k}\binc{n+1-k}{j}
\lambda^{k}d^{n+1-j}(x)d^{k+j}(y)\\
&& + \sum_{j=0}^{n+1}\binc{n+1}{j}d^{n+1-j}(x)d^{j}(y)
+ \lambda^{n+1}d^{n+1}(x)d^{n+1}(y)\\
&=& \sum_{k=0}^{n+1}\sum_{j=0}^{n+1-k}\binc{n+1}{k}\binc{n+1-k}{j}
\lambda^{k}d^{n+1-j}(x)d^{k+j}(y).
\end{eqnarray*}
}
This completes the induction.
\end{proof}
\medskip

We now briefly study degenerated $\lambda$-differential operators, that is, weak differential operators $d$ for which $d(1)\neq 0$.

We first note that, for any $\lambda\in \bfk$ and any
$\bfk$-algebra $R$, the zero map
$$d: R\to R,\ d(r)=0,\ \forall\, r\in R$$
is a differential operator of weight $\lambda$,
called the {\bf zero differential operator}
of weight $\lambda$.

We next note that for any $\lambda\in \bfk$ that is invertible and for any $\bfk$-algebra $R$, the map
\begin{equation}
d: R\to R,\ d(r)=-\lambda^{-1} r,\ \forall\, r\in R,
\mlabel{eq:scadiff}
\end{equation}
is a weak differential operator of weight $\lambda$.
We call such an operator (resp. algebra)
a {\bf scalar differential operator (resp. algebra) of weight $\lambda$}.
We remark that by our definition, the zero map is not a scalar differential operator even though the zero map is given by a scalar multiplication.

For $\lambda\in \bfk$ invertible, it is also easy to check that
$$ P_\lambda: R\to R, P_\lambda(r)=-\lambda\, r,\quad \forall\, r\in\, R,$$
is a Rota-Baxter operator of weight $\lambda$. Further
$d_\lambda \circ P_\lambda=\id$. This gives an instance of a degenerated differential Rota-Baxter algebra of weight $\lambda$.

\begin{prop}
Let $\lambda \in \bfk$. Let $(R,d)$ be a weak differential $\bfk$-algebra of weight $\lambda$ with no zero divisors. Then the following statements are equivalent.
\begin{enumerate}
\item
$\lambda$ is invertible and $d$ is a scalar differential operator of weight $\lambda$. \mlabel{it:sca}
\item
$\lambda$ is invertible and $d(1)=-\lambda^{-1}.$
\mlabel{it:1sca}
\item
$d(1)\neq 0$.
\mlabel{it:dnonzero}
\item
For every $r\in R$, $d(r)$ is a non-zero $\bfk$-multiple of $r$.
\mlabel{it:dep}
\end{enumerate}
\mlabel{pp:difftype}
\end{prop}
\begin{proof}
We clearly have
(\mref{it:sca}) $\Rightarrow$ (\mref{it:1sca})
$\Rightarrow$ (\mref{it:dnonzero})
and
(\mref{it:sca}) $\Rightarrow$ (\mref{it:dep}).
So we only need to prove
(\mref{it:dnonzero}) $\Rightarrow$ (\mref{it:sca})
and (\mref{it:dep}) $\Rightarrow$ (\mref{it:dnonzero}).

\noindent
{\bf (\mref{it:dnonzero}) $\Rightarrow$ (\mref{it:sca}): }
By Eq.~(\mref{eq:diff}), for any $x\in R$, we have
$$ d(x)=d(1)x+d(x)+\lambda d(1)d(x).$$
Thus $d(1)\big(x+\lambda d(x)\big)=0.$
Since $R$ has no zero divisors, if $d(1)\neq 0$, then we have
\begin{equation}
x+\lambda d(x)=0.
\mlabel{eq:diffsc}
\end{equation}
Letting $x=1$, we have
$ 0=1+\lambda d(1) = 1+ d(1) \lambda$ since $\lambda\in \bfk$. Thus $\lambda$ is invertible and Eq.~(\mref{eq:diffsc}) gives $d(x)=-\lambda^{-1} x,$ as needed.

\noindent
{\bf (\mref{it:dep}) $\Rightarrow$ (\mref{it:dnonzero}): }
Taking $r=1$ in (\mref{it:dep}) we have
$d(1)=\alpha(1)\neq 0$. Thus we have (\mref{it:dnonzero}).
\end{proof}

\begin{coro}
Let $(R,d)$ be a weak differential algebra of weight $\lambda$ that has no zero divisors. If $d$ is not a scalar differential operator, then it is a differential operator of weight $\lambda$.
\mlabel{co:nonsc}
\end{coro}

\begin{remark}{\rm
Since a scalar differential algebra is just an algebra with a fixed scalar multiplication, its study can be reduced to the study of algebras.
By Corollary~\mref{co:nonsc}, a non-scalar weak differential algebra is a differential algebra under a mild restriction. This justifies the requirement in our Definition~\mref{de:main} that a $\lambda$-differential algebra be nondegenerated.
A more careful study of degenerated differential operators will be carried out in another study.
}
\mlabel{re:deg}
\end{remark}

\subsection{Free differential algebras of weight $\lambda$}
Using the same construction as for free differential algebras (of weight 0), we obtain free differential algebras of weight $\lambda$ in both the commutative and non-commutative case.

\begin{theorem}
Let $X$ be a set. Let
$$\diffs(X)=X\times \NN= \{ x^{(n)}\, \big|\, x\in X, n\geq 0\}.$$
\begin{enumerate}
\item
Let $\bfk\diffa{X}$ be the free commutative algebra $\bfk[\diffs{X}]$ on the set $\diffs{X}$.
Define $d_X: \bfk\diffa{X} \to \bfk\diffa{X}$ as follows. Let $w=u_1\cdots u_k, u_i\in \diffs{X}$,
$1\leq i\leq k$, be a commutative word from the alphabet set $\Delta(X)$.  If $k=1$, so that $w=x^{(n)}\in \Delta(X)$, define $d_X(w)=x^{(n+1)}$. If $k>1$, recursively define
\begin{equation}
 d_X(w)=d_X(u_1)u_2\cdots u_k+u_1d_X(u_2\cdots u_k)+\lambda d_X(u_1)d_X(u_2\cdots u_k).
 \mlabel{eq:prodind}
 \end{equation}
Further define $d_X(1)=0$ and then extend $d_X$ to $\bfk\diffa{X}$ by linearity.
Then $(\bfk\diffa{X},d_X)$ is the free commutative differential algebra of weight $\lambda$ on the set $X$.
\mlabel{it:commfreediff}
\item
Let $\bfk^{NC}\diffa{X}$ be the free noncommutative algebra $\bfk^{NC}[\diffs{X}]$ on the set $\diffs{X}$.
Define $d^{NC}_X: \bfk^{NC}\diffa{X} \to \bfk^{NC}\diffa{X}$ on the noncommutative words from the alphabet set $\diffs{X}$ in the same way as $d_X$ is defined in (\mref{it:commfreediff}).
Then $(\bfk^{NC}\diffa{X},d_X^{NC})$ is the free noncommutative differential algebra of weight $\lambda$ on the set $X$.
\mlabel{it:ncfreediff}
\end{enumerate}
\mlabel{thm:freediff}
\end{theorem}
\begin{proof}
We just give a proof of (\mref{it:commfreediff}). The proof of (\mref{it:ncfreediff}) is the same. In either case, it is similar to the proof of the $\lambda=0$ case~\cite[p.70]{Kol}.

Let $(R,d)$ be a commutative $\lambda$-differential algebra and let $f:X\to R$ be a set map. We extend $f$ to a $\lambda$-differential algebra homomorphism $\free{f}: \bfk\diffa{X} \to R$ as follows.

Let $w=u_1\cdots u_k, u_i\in \diffs{X}$,
$1\leq i\leq k$, be a commutative word from the alphabet set $\diffs{X}$.  If $k=1$, then $w=x^{(n)}\in \diffs{X}$. Define
\begin{equation}
\free{f}(w)=d^n(f(x)).
\mlabel{eq:diffred1}
\end{equation}
Note that this is the only possible definition in order for $\free{f}$ to be a $\lambda$-differential algebra homomorphism. If $k>1$, recursively define
$$ \free{f}(w)=\free{f}(u_1)\free{f}(u_2\cdots u_k).$$
Further define $\free{f}(1)=1$ and then extend $\free{f}$ to $\bfk\diffa{X}$ by linearity. This is the only possible definition in order for $\free{f}$ to be an algebra homomorphism.

Since $\bfk\diffa{X}$ is the free commutative algebra on
$\diffs{X}$, $\free{f}$ is an algebra homomorphism.
So it remains to verify that, for all commutative words $w=u_1\cdots u_k$ from the alphabet set $\diffs{X}$,
\begin{equation}
 \free{f} (d_X(w))=d(\free{f}(w)),
\mlabel{eq:diffind1}
\end{equation}
for which we use induction on $k$. The case when $k=1$ follows immediately from Eq.~(\mref{eq:diffred1}). For the inductive step,  by Eq.~(\mref{eq:prodind}):
\begin{eqnarray*}
\free{f}( d_X(w))&=&\free{f}(d_X(u_1)u_2\cdots u_k) + \free{f}(u_1d_X(u_2\cdots u_k))+\lambda \free{f}(d_X(u_1)d_X(u_2\cdots u_k))\\
&=& \free{f}(d_X(u_1))\free{f}(u_2\cdots u_k) + \free{f}(u_1)\free{f}(d_X(u_2\cdots u_k))+\lambda \free{f}(d_X(u_1))\free{f}(d_X(u_2\cdots u_k)).
\end{eqnarray*}
Then by Eq.~(\mref{eq:diffred1}), the induction hypothesis on $k$ and the $\lambda$-differential algebra relation for $d$, the last sum equals to
$d(\free{f}(w))$.
\end{proof}

\subsection{Cofree differential algebras of weight $\lambda$}

For any $\bfk$-algebra $A$, let $A^{\NN}$ denote the $\bfk$-module of all
functions $f:\NN \rightarrow A$.  We define a product on $A^{\NN}$
by defining, for any $f, g \in A^{\NN}$, $fg \in A^{\NN}$ by
\[(fg)(n) = \sum_{k=0}^{n}\sum_{j=0}^{n-k}\binc{n}{k}\binc{n-k}{j}
\lambda^{k}f(n-j)g(k+j).\]
Note that this definition is motivated by Proposition~\mref{pp:basic}.(\mref{it:prod}).
It is easily checked that this product is commutative, associative,
distributive over addition, and has an identity $\bfone_{A^{\NN}}$
defined by $\bfone_{A^{\NN}}(n) = 0$ if $n \neq 0$ and
$\bfone_{A^{\NN}}(0) = \bfone_{A}$.  We call this product the
{\bf $\lambda$-Hurwitz product} on $A^{\NN}$, since if we take
$\lambda = 0$, the product reduces to
\[(fg)(n) = \sum_{k = 0}^{n}\binc{n}{k}f(n-k)g(k),\]
which is the usual Hurwitz product defined in~\cite{Ke}.
We denote the $\bfk$-algebra
$A^{\NN}$ with this product by $DA$, and call it the $\bfk$-algebra of
{\bf $\lambda$-Hurwitz series over $A$}.
Also, there is, for any
$\bfk$-algebra $A$, a homomorphism $\kappa_{A}:A \rightarrow DA$
of $\bfk$-algebras defined by $\kappa_A(a) = a\bfone_{A^{\NN}}$.
This makes $DA$ into an $A$-algebra, where for any $a \in A$ and
any $f \in DA$, $af \in DA$ is given by $(af)(n) = a(f(n))$.

The $\bfk$-algebra $DA$ behaves much like the ring of Hurwitz series.
The following proposition is one instance of this.  We first define
a map
\begin{equation}
\partial_A:DA \rightarrow DA, \quad (\partial_A(f))(n) = f(n+1), n\in \NN, f\in DA.
\mlabel{eq:hurw}
\end{equation}
\begin{prop}
\mlabel{prop:deriv}
The map $\partial_A$ is a $\lambda$-derivation on $DA$.
\end{prop}
\proof
It is clear that $\partial_A$ is a mapping of $\bfk$-modules, so
all that remains is to show that for any $f,g \in DA$,
\[\partial_A(fg) = \partial_A(f)g + f\partial_A(g)
+ \lambda\partial_A(f)\partial_A(g).\]
But because of the definition
of the $\lambda$-Hurwitz product, the proof of this equation is
virtually identical to the proof of Proposition~\ref{pp:basic}.(\mref{it:prod}) and
is left to the reader.
\proofend

It follows from Proposition~\ref{prop:deriv} that $(DA, \partial_A)$
is a $\lambda$-differential $\bfk$-algebra.  If $h:A \rightarrow B$ is a
$\bfk$-algebra homomorphism, one checks that $Dh:DA \rightarrow DB$ defined
by $((Dh)(f))(n) = h(f(n))$ is a morphism of $\bfk$-algebras, and that
$\partial_B \circ Dh = Dh \circ \partial_A$.  Recalling
that $\Dif=\Dif_\lambda$ denotes the category
of $\lambda$-differential $\bfk$-algebras, we see that we have a functor
$G:\Alg_{\bfk} \rightarrow \Dif$ given on objects
$A \in \Alg $ by $G(A) = (DA, \partial_A)$ and on morphisms
$h:A \rightarrow B$ in $\Alg $ by $G(h) = Dh$ as defined
above.  Letting $V:\Dif \rightarrow \Alg $ denote the forgetful
functor defined on objects $(R,d) \in \Dif$ by $V(R,d) = R$ and
on morphisms $f:(R,d) \rightarrow (S,e)$ in $\Dif$ by $V(f) = f$, we
have the following characterization of $G(A)$.

\begin{prop}
\mlabel{prop:cofree}
The functor $G:\Alg  \rightarrow \Dif$ defined above is the right
adjoint of the forgetful functor $V:\Dif \rightarrow \Alg $.
\end{prop}
\proof
By~\cite{Ma}, it is equivalent to show that there are two natural
transformations $\eta: id_{\Dif} \rightarrow GV$ and
$\varepsilon: VG \rightarrow id_{\Alg }$ satisfying the equations
$G\varepsilon \circ \eta G = G$ and $\varepsilon V \circ V\eta = V$.
Here $id_{\Dif}$ denotes the
identity functor on $\Dif$, and similarly for $id_{\Alg }$.

For any $A \in \Alg $, define $\varepsilon_A: DA \rightarrow A$
for any $f \in DA$ by $\varepsilon_A(f) = f(0).$  One checks that
$\varepsilon_A$ is a morphism of $\bfk$-algebras, and that if
$h:A \rightarrow B$ is any morphism of $\bfk$-algebras, then
$\varepsilon_B \circ Dh = h \circ \varepsilon_A$, i.e., $\varepsilon$
is a natural transformation as desired.

For any $(R,d) \in \Dif$, $x \in R$ and $n \in \NN$, define
$\eta_{(R,d)}:(R,d) \rightarrow (DR, \partial_R)$ by
$(\eta_{(R,d)}(x))(n) = d^{(n)}(x)$.  It is not difficult to see
that $\eta_{(R,d)}$ is $\bfk$-linear, and it is immediate from
Proposition~\mref{pp:basic}.(\mref{it:prod}) that for any $x, y \in R$,
$(\eta_{(R,d)}(x))(\eta_{(R,d)}(y)) = (\eta_{(R,d)}(xy)).$
Also, it is clear that $\partial_R \circ \eta_{(R.d)} =
\eta_{(R.d)} \circ d$, so that $\eta_{(R,d)}$ is a morphism in
$\Dif$.
Further, if $f:(R,d) \rightarrow (S,e)$ is a morphism in $\Dif$,
then one sees that $\eta_{(S,e)} \circ f = Df \circ \eta_{(R,d)}.$
Hence $\eta$ is a natural transformation.

To see that $G\varepsilon \circ \eta G = G$, let $A \in \Alg $,
$f \in DA$ and $n \in \NN$.  Then
$$((D\varepsilon_A)(\eta_{(DA,\partial_A)}(f)))(n) =
\varepsilon_A(\eta_{(DA,\partial_A)}(f)(n)) =
\varepsilon_A(\partial_A^{(n)}(f)) =
(\partial_A^{(n)}(f))(0) =f(n).$$

Similarly, to see that $\varepsilon V \circ V\eta = V$,
let $(R,d) \in \Dif$, and $x \in R$.  Then
$\varepsilon_R(\eta_{(R,d)}(x)) =
(\eta_{(R,d)}(x))(0) = d^{(0)}(x) = x.$
\proofend

The following corollary gives a ``universal mapping property''
characterization of the $\lambda$-differential $\bfk$-algebra of
$\lambda$-Hurwitz series as the cofree $\lambda$-differential
$\bfk$-algebra on any $\bfk$-algebra $A$.

\begin{coro}
\mlabel{coro:cofree}
Let $(R,d)$ be any $\lambda$-differential $\bfk$-algebra, and let $A$
be any $\bfk$-algebra.  For any $\bfk$-algebra homomorphism
$f:R \rightarrow A$, there is a unique morphism of $\lambda$-differential $\bfk$-algebras $\tilde{f}:(R,d) \rightarrow (DA,\partial_A)$
such that $\varepsilon_A \circ V\tilde{f} = f.$
\end{coro}
\begin{proof}
This follows from page 81, Theorem 2 in~\cite{Ma}.
\end{proof}

\delete{
Define $\tilde{f}$ to be the composite $Gf \circ \eta_{(R,d)}$.
Then $\varepsilon_A \circ V\tilde{f} =
\varepsilon_A \circ VGf \circ V\eta_{(R,d)} =
f \circ \varepsilon_{V(R,d)} \circ V\eta_{(R,d)} = f$.
Suppose further that $g: (R,d) \rightarrow (DA,\partial_A)$ is
any morphism of $\lambda$-differential $\bfk$-algebras such that
$\varepsilon_A \circ Vg = f$.  Then $Gf \circ \eta_{(R,d)}
= G\varepsilon_A \circ GVg \circ \eta_{(R,d)}
= G\varepsilon_A \circ \eta_{GA} \circ g = g$.
\proofend
}

We next show that $DA$ provides another example of differential Rota-Baxter algebras.  Define
\begin{equation}
\pi_A: DA\to DA, \quad (\pi_A(f))(n)=f(n-1), n\geq 1,
(\pi_A(f))(0)=0, f\in DA.
\mlabel{eq:hurwrb}
\end{equation}
\begin{prop}
The triple $(DA,\partial_A, \pi_A)$ is a differential Rota-Baxter algebra of weight $\lambda$.
\mlabel{pp:hurwdrb}
\end{prop}
\begin{proof}
Since
$$(\pi_A (\partial_A(f)))(n)=(\pi_A(f))(n+1)=f(n)$$
for $f\in DA$, we have $\pi_A\circ \partial_A=\id_{DA}$.
Thus we only need to verify that $\pi_A$ is a Rota-Baxter operator of weight $\lambda$. Let $H\in DA$ be defined by
\begin{equation}
H=\pi_A(f)\pi_A(g)-\pi_A(\pi_A(f)g)-\pi_A(f\pi_A(g))- \lambda \pi_A(fg).
\mlabel{eq:pirb}
\end{equation}
By Proposition~\mref{prop:deriv}, we have
$\partial_A(H)=0$. Thus $H$ is of the form $H:\NN\to A$ with $H(n)=0, n>0$ and $H(0)=k$ for some $k\in \bfk$.
But by definition, $\pi_A(0)=0$. Thus $H(0)=0$ and so $H=0$. This shows that $\pi_A$ is a Rota-Baxter operator of weight $\lambda$.
\end{proof}

\section{Free commutative differential Rota-Baxter algebras}
\mlabel{sec:commfree}
We briefly recall the construction of free commutative Rota-Baxter algebras in terms of mixable shuffles~\mcite{G-K1,G-K2}. Let $A$ be a commutative $\bfk$-algebra. Define
\[ \sha (A)= \bigoplus_{k\in\NN}
    A^{\otimes (k+1)}
= A\oplus A^{\otimes 2}\oplus \ldots. \]
Let $\fraka=a_0\ot \cdots \ot a_m\in A^{\ot (m+1)}$ and
$\frakb=b_0\ot \cdots \ot b_n\in A^{\ot (n+1)}$. If $m=0$ or $n=0$, define
\begin{equation}
\fraka \shpr \frakb =\left \{\begin{array}{ll}
    (a_0b_0)\ot b_1\ot \cdots \ot b_n, & m=0, n>0,\\
    (a_0b_0)\ot a_1\ot \cdots \ot a_m, & m>0, n=0,\\
    a_0b_0, & m=n=0.
    \end{array} \right .
\end{equation}
If $m>0$ and $n>0$, recursively define
\begin{eqnarray}
\fraka \shpr \frakb & = &
(a_0b_0)\ot \Big(
a_1\ot \big( (a_2\ot \cdots \ot a_m) \shpr (b_1\ot \cdots \ot b_n)\big) \notag \\
&&
\qquad \qquad +
\;b_1\ot \big( (a_1\ot \cdots \ot a_m) \shpr (b_2\ot \cdots \ot b_n)\big) \mlabel{eq:shpr}\\
&& \qquad \qquad +
\lambda\, a_1b_1\ot \big( (a_2\ot \cdots \ot a_m) \shpr (b_2\ot \cdots \ot b_n)\big)\Big)
\notag
\end{eqnarray}
with the convention that
\begin{eqnarray*}
&&(a_2\ot \cdots \ot a_m) \shpr (b_1\ot \cdots \ot b_n)
=b_1\ot \cdots \ot b_n, {\rm\ if\ } m=1, n>1;\\
&&(a_1\ot \cdots \ot a_m) \shpr (b_2\ot \cdots \ot b_n)
=a_1\ot \cdots \ot m_n, {\rm\ if\ } m>1, n=1;\\
&&a_1b_1\ot \big( (a_2\ot \cdots \ot a_m) \shpr (b_2\ot \cdots \ot b_n)\big)=a_1b_1, {\rm\ if\ } m=n=1.
\end{eqnarray*}
Extending by additivity, we obtain a $\bfk$-bilinear map

\[ \shpr: \sha (A) \times \sha (A) \rar \sha (A), \]
called the mixable shuffle product on $\sha(A)$.
Define a $\bfk$-linear endomorphism $P_A$ on
$\sha (A)$ by assigning
\[ P_A( x_0\otimes x_1\otimes \ldots \otimes x_n)
=\bfone_A\otimes x_0\otimes x_1\otimes \ldots\otimes x_n, \]
for all
$x_0\otimes x_1\otimes \ldots\otimes x_n\in A^{\otimes (n+1)}$
and extending by additivity.
Let $j_A:A\rar \sha (A)$ be the canonical inclusion map.

\begin{theorem}
\mlabel{thm:shua}
The pair $(\sha (A),P_A)$, together with the natural embedding
$j_A:A\rightarrow \sha (A)$,
is a free commutative Rota-Baxter $\bfk$-algebra on $A$ of weight $\lambda$.
In other words,
for any Rota-Baxter $\bfk$-algebra $(R,P)$ and any
$\bfk$-algebra map
$\varphi:A\rar R$, there exists
a unique Rota-Baxter $\bfk$-algebra homomorphism
$\tilde{\varphi}:(\sha (A),P_A)\rar (R,P)$ such that
$\varphi = \tilde{\varphi} \circ j_A$ as $\bfk$-algebra homomorphisms.
\end{theorem}
It is proved in~\cite{G-K1} that $\sha(A)$ with the mixable shuffle product is the free commutative Rota-Baxter algebra on $A$. The mixable shuffle product is shown to be the same as the quasi-shuffle product~\mcite{E-G1,G-Z,Ho}.

Since $\shpr$ is compatible with the multiplication in $A$, we will often suppress
the symbol $\shpr$ and simply denote $x y$ for $x\shpr y$
in $\sha (A)$, unless there is a danger of confusion.

Let $(A,d_0)$ be a commutative differential $\bfk$-algebra of weight $\lambda$.
Define an operator $d_A$ on $\sha (A)$ by assigning
\begin{eqnarray*}
\lefteqn{ d_A(x_0\otimes x_1\otimes\ldots\otimes x_n)}\\
&=&
d_0(x_0)\otimes x_1\otimes \ldots \otimes x_n + x_0x_1\otimes
x_2 \otimes \ldots \otimes x_n +\lambda d_0(x_0) x_1\otimes x_2\otimes
\ldots \otimes x_n
\end{eqnarray*}
for $x_0\otimes \ldots \otimes x_n\in
A^{\otimes (n+1)}$ and then extending by $\bfk$-linearity. Here we
use the convention that when $n=0$, $d_A(x_0)=d_0(x_0)$.

\begin{theorem}
\mlabel{thm:commdiffrb}
Let $(A,d_0)$ be a commutative
differential $\bfk$-algebra of weight $\lambda$.
\begin{enumerate}
\item
$(\sha (A),d_A,P_A)$ is a commutative differential Rota-Baxter $\bfk$-algebra of weight $\lambda$. The $\bfk$-algebra embedding
\[ j_A: A \rar \sha (A)\]
is a morphism of differential $\bfk$-algebras of weight $\lambda$.
\mlabel{it:freecommdiffa}
\item
The quadruple $(\sha (A),d_A,P_A,j_A)$ is a free commutative differential
Baxter $\bfk$-algebra of weight $\lambda$ on $(A,d_0)$, as described by the
following universal property: For any commutative differential Rota-Baxter
$\bfk$-algebra $(R,d,P)$ of weight $\lambda$ and any $\lambda$-differential $\bfk$-algebra map $\varphi: (A,d_0)\rar (R,d)$, there
exists a unique $\lambda$-differential Rota-Baxter $\bfk$-algebra
homomorphism $\tilde{\varphi}:(\sha (A),d_A,P_A)\rar (R,d,P)$
such that the diagram
\[\xymatrix{
(A,d_0) \ar[r]^-{j_A} \ar[dr]_{\varphi}
    & (\sha (A),d_A) \ar[d]^{\tilde{\varphi}} \\
& (R,d) } \] commutes in the category of commutative differential
$\bfk$-algebras of weight $\lambda$.
\mlabel{it:freecommdiffrba}
\item
Let $X$ be a set and let $\bfk\diffa{X}$ be the free commutative differential algebra of weight $\lambda$ on $X$. The quadruple $(\sha (\bfk\diffa{X}),d_{\bfk\diffa{X}},P_{\bfk\diffa{X}},j_X)$ is a free commutative differential
Baxter $\bfk$-algebra of weight $\lambda$ on $X$, as described by the following universal property: For any commutative differential Rota-Baxter
$\bfk$-algebra $(R,d,P)$ of weight $\lambda$ and any set map $\varphi: X\to R$, there
exists a unique $\lambda$-differential Rota-Baxter $\bfk$-algebra
homomorphism $\tilde{\varphi}:(\sha (\bfk\diffa{X}),d_{\bfk\diffa{X}},
P_{\bfk\diffa{X}})\rar (R,d,P)$
such that $\tilde{\varphi}\circ j_X = \varphi.$
\mlabel{it:freecommdiffrbx}
\end{enumerate}
\end{theorem}

\begin{proof} (\mref{it:freecommdiffa}).
 For any $x=x_0\otimes \ldots\otimes x_m\in A^{\otimes
(m+1)}$, by definition we have
$$
d_A (P_A(x))= d_A(1\otimes x_0\otimes \ldots \otimes
x_m)= x_0\otimes \ldots \ot x_m.
$$
Thus $d_A\circ P_A$ is the identity map on $\sha(A)$.
So it remains to prove that for any $m,\ n\in
\NN_+$ and $x=x_0\otimes \ldots \otimes x_m\in A^{\otimes (m+1)}$
and $y=y_0\otimes \ldots \otimes y_n\in A^{\otimes (n+1)}$, we
have
\begin{equation}
 d_A(x\shpr y) = d_A(x) \shpr y+ x\shpr d_A(y)
+\lambda d_A(x)\shpr d_A(y).
\mlabel{eq:dif}
\end{equation}
If $m=0$ or $n=0$, then the equation follows from the definition
of $d_A$. Now consider the case when $m,\ n\in \NN_+$. Denoting
$x^+= x_1\otimes \ldots \otimes x_{m}$ and $y^+= y_1 \otimes
\ldots \otimes y_{m}$, we have $x= x_0\shpr P_A(x^+)$, $y=
y_0\shpr P_A(y^+)$ and Eq.~(\mref{eq:shpr}) can be rewritten as
\begin{eqnarray*}
 x\shpr y &=& (x_0 y_0)\shpr (P_A(x^+) \shpr P_A(y^+)) \\
 &=& (x_0 y_0)\shpr \left ( P_A(x^+\shpr P_A(y^+))+P_A(y^+\shpr P_A(x^+))
    + \lambda P_A(x^+ \shpr y^+)\right ) \\
 &=& (x_0 y_0) \shpr P_A(x^+\shpr P_A(y^+)+y^+\shpr P_A(x^+)
    +\lambda (x^+\shpr y^+)).
\end{eqnarray*}
It follows from the definition of $d_A$ that, for any $z_0\in A$
and $z\in \sha (A)$,
\begin{eqnarray*}
d_A(z_0\shpr P_A(z))&=& d_0(z_0)\shpr P_A(z) + z_0\shpr
d_A(P_A(z))
    +\lambda d_0(z_0)\shpr d_A(P_A(z))\\
&=& d_0(z_0) \shpr P_A(z) + z_0 \shpr z
    +\lambda d_0(z_0) \shpr z.
\end{eqnarray*}
Then
\begin{eqnarray*}
\lefteqn{d_A  (x\shpr y)}\\
 &=&
    d_A( (x_0 y_0)\shpr  P_A(x^+\shpr P_A(y^+)
    +y^+\shpr P_A(x^+)
    + \lambda (x^+ \shpr y^+)))\\
 &=& d_0(x_0y_0) \shpr P_A(x^+\shpr P_A(y^+)
    +y^+\shpr P_A(x^+)
    + \lambda (x^+ \shpr y^+))\\
 &&+ (x_0 y_0) \shpr (x^+\shpr P_A(y^+)
    +y^+\shpr P_A(x^+)
    + \lambda (x^+ \shpr y^+) )\\
 &&+ \lambda d_0(x_0 y_0) \shpr (x^+\shpr P_A(y^+)
    +y^+\shpr P_A(x^+)
    + \lambda (x^+ \shpr y^+))\\
&=&(d_0(x_0)y_0 +x_0 d_0(y_0)+\lambda d_0(x_0)d_0(y_0)) \shpr
(P_A(x^+) \shpr P_A(y^+))\\ &&+(x_0 y_0) \shpr (x^+\shpr
P_A(y^+)+y^+\shpr P_A(x^+)
    +x^+ \shpr y^+)\\
&& +\lambda (d_0(x_0)y_0 +x_0 d_0(y_0)+\lambda d_0(x_0)
d_0(y_0))\\ &&  \shpr (x^+\shpr P_A(y^+)+y^+\shpr P_A(x^+)+x^+
\shpr y^+).
\end{eqnarray*}
Also
\begin{eqnarray*}
\lefteqn{ x\shpr d_A(y) +y\shpr d_A(x)+\lambda d_A(x)\shpr
d_A(y)}\\ &=& (x_0 \shpr P_A(x^+)) \shpr d_A(y_0\shpr P_A(y^+)) +
(y_0\shpr P_A(y^+))\shpr d_A(x_0\shpr P_A(x^+))\\ && + \lambda
d_A(x_0\shpr P_A(x^+))\shpr d_A(y_0\shpr P_A(y^+))\\ &=& (x_0\shpr
P_A(x^+))\shpr \left (d_0(y_0) \shpr P_A(y^+) + y_0\shpr y^+
+\lambda d_0(y_0)\shpr y^+)
    \right )\\
&&+ (y_0\shpr P_A(y^+))\shpr \left (d_0(x_0) \shpr P_A(x^+) +
x_0\shpr x^+ +\lambda d_0(x_0)\shpr x^+
    \right )\\
&&+\lambda \left (d_0(x_0) \shpr P_A(x^+) x_0\shpr x^+ +\lambda
d_0(x_0)\shpr x^+
    \right )\\
&&\ \ \ \shpr \left (d_0(y_0) \shpr P_A(y^+) y_0\shpr y^+ +\lambda
d_0(y_0)\shpr y^+
    \right ).
\end{eqnarray*}
Comparing last terms of the above two equations, we see that
equation~(\ref{eq:dif}) holds.

The second statement follows directly from the
definition of $d_A$.
\medskip

\noindent
(\mref{it:freecommdiffrba}).
Now let $(R,d,P)$ be a commutative differential Rota-Baxter
$\bfk$-algebra of weight $\lambda$ and let $\varphi: (A,d_0)\rar (R,d)$ be a $\lambda$-differential
$\bfk$-algebra map. Then in particular $\varphi$ is a $\bfk$-algebra
map. So by Theorem~\ref{thm:shua}, there is a unique Rota-Baxter
$\bfk$-algebra map $\tilde{\varphi}: (\sha (A),P_A)\rar (R,P)$ such
that
\begin{equation}
 \varphi=\tilde{\varphi}\circ j_A
\end{equation}
in the category of $\bfk$-algebras. We next show that
$\tilde{\varphi}$ is a differential $\bfk$-algebra map.

For any $x_0\otimes x_1\otimes \ldots \otimes x_n\in A^{\otimes
(n+1)}$, we have
\begin{eqnarray*}
\lefteqn{ d (\tilde{\varphi}(x_0\otimes x_1\otimes \ldots
    \otimes x_n))
= d(\tilde{\varphi}(x_0 \shpr P_A(x_1\otimes\ldots\otimes
x_n)))}
\\ &=& d(\tilde{\varphi}(x_0)
    \tilde{\varphi}(P_A(x_1\otimes\ldots \otimes x_n)))\\
&=& d(\tilde{\varphi}(x_0)
    P(\tilde{\varphi}(x_1\otimes\ldots\otimes x_n)))\\
&=& d(\tilde{\varphi}(x_0))
    P(\tilde{\varphi}(x_1\otimes\ldots\otimes x_n))
+ \tilde{\varphi}(x_0)
    d(P(\tilde{\varphi}(x_1\otimes\ldots\otimes x_n)))\\
&& + \lambda d(\tilde{\varphi}(x_0))
    d(P(\tilde{\varphi}(x_1\otimes\ldots\otimes x_n)))\\
&=& d(\tilde{\varphi}(x_0))
    \tilde{\varphi}(P_A(x_1\otimes\ldots\otimes x_n))
+ \tilde{\varphi}(x_0)
    \tilde{\varphi}(x_1\otimes\ldots\otimes x_n)\\
&& + \lambda d(\tilde{\varphi}(x_0))
    \tilde{\varphi}(x_1\otimes\ldots\otimes x_n).
\end{eqnarray*}
Also
\begin{eqnarray*}
\lefteqn{\tilde{\varphi}(d_A(x_0\otimes x_1\otimes \ldots\otimes x_n)) = \tilde{\varphi}(d_A(x_0\shpr
P_A(x_1\otimes\ldots\otimes x_n)))}\\
&=&
\tilde{\varphi}(d_A(x_0)\shpr P_A(x_1\otimes\ldots\otimes x_n) +
x_0 \shpr d_A(P_A(x_1\otimes\ldots\otimes x_n))\\
&&+ \lambda\,d_A(x_0)
\shpr d_A(P_A(x_1\otimes\ldots\otimes x_n)))\\
&=&
\tilde{\varphi}(d_A(x_0)\shpr P_A(x_1\otimes\ldots\otimes x_n) +
x_0 \shpr (x_1\otimes\ldots\otimes x_n)\\
&&+ \lambda\, d_A(x_0) \shpr
(x_1\otimes\ldots\otimes x_n))\\
&=& \tilde{\varphi}(d_A(x_0))
    \tilde{\varphi}(P_A(x_1\otimes\ldots\otimes x_n))
+ \tilde{\varphi}(x_0)\tilde{\varphi}(x_1\otimes\ldots\otimes
x_n)\\
&&
+\lambda\,\tilde{\varphi}(d_A(x_0))\tilde{\varphi}(x_1\otimes\ldots\otimes
x_n).
\end{eqnarray*}
Since
\begin{eqnarray*}
d(\tilde{\varphi}(x_0))&=&d(\tilde{\varphi}(j_A(x_0)))
= d(\varphi(x_0))= \varphi(d_0(x_0)) \\
&=&
\tilde{\varphi}(j_A(d_0(x_0))) =
\tilde{\varphi}(d_A(j_A(x_0))) = \tilde{\varphi}(d_A(x_0)),
\end{eqnarray*}
we have proved that
\[d (\tilde{\varphi}(x_0\otimes x_1\otimes \ldots
    \otimes x_n))=
\tilde{\varphi}(d_A(x_0\otimes x_1\otimes \ldots\otimes x_n)).\]
This shows that $\tilde{\varphi}$ is a differential $\bfk$-algebra
homomorphism. Since $\varphi$ and $j_A$ are differential
$\bfk$-algebra homomorphisms, we see that equation~(\ref{eq:diff})
holds in the category of differential $\bfk$-algebras.
\smallskip

\noindent
(\mref{it:freecommdiffrbx}). The forgetful functor from the category $\DRB_\lambda$ to the category $\Set$ of sets is the composition of the forgetful functors from $\DRB_\lambda$ to $\Alg$ and from $\Alg$ to $\Set$. By Theorem 1 in page 101 of~\mcite{Ma}, the adjoint functor of a composed functor is the composition of the adjoint functors. This proves (\mref{it:freecommdiffrbx}).
\end{proof}

\delete{
For a given $(R,d,P)$ in $\DRB_\lambda$ and a set
map $\varphi:X\rar R$,
by the universal property of the free $\lambda$-differential algebra
$(\bfk\diffa{X},d_X)$, $\varphi$ extends uniquely to a
$\lambda$-differential algebra homomorphism
$\bar{\varphi}:(\bfk\diffa{X},d_X) \rar (R,d)$.
Then by item (\mref{it:freecommdiffrba}), $\bar{\varphi}$
extends uniquely to a homomorphism of $\lambda$-differential Rota-Baxter $\bfk$-algebras
\[\tilde{\varphi}: (\sha \diffa{X},d_{\bfk\diffa{X}},P_{\bfk\diffa{X}})\to (R,d,P).\]
This proves the proposition.
}

\section{Free noncommutative differential Rota-Baxter algebras}
\mlabel{sec:ncfree}
We now consider the noncommutative analog of Section~\mref{sec:commfree}.

\subsection{Free noncommutative Rota-Baxter algebras}
We first summarize the construction of free noncommutative Rota-Baxter algebras on a set $X$ in terms of angularly decorated planar rooted trees. See~\mcite{E-G0} (as well as~\mcite{A-M}) for further details.

\subsubsection{Rota-Baxter algebra on rooted trees}
\mlabel{sss:rbtree}
We follow the notations and terminologies in~\cite{Di,We}. A free tree is an undirected graph that is connected and
contains no cycles. A {\bf rooted tree} is a free tree in which a particular vertex has been distinguished as the {\bf root}. A {\bf planar rooted tree} is a rooted tree with a fixed embedding into the plane.
For example,
$$
 \ta1 \;\quad
 \tb2   \quad\;
 \tc3   \quad\;
 \td31  \quad\;
 \te4   \quad\;
 \tf41  \quad\;
 \tg42 \;\quad
 \thj44 \quad\;
 \th43 \;\quad
\delete{
  \ti5  \quad\;
 \tj51  \quad\;
 \tk52  \quad\;
 \tl53  \quad\;
 \tm54  \quad\;
 \tn55
} \cdots
$$
The depth $\depth(T)$ of a rooted tree $T$ is the length of the longest path from its root to its leaves.

Let $\calt$ be the set of planar rooted trees.
A {\bf planar rooted forest} is a noncommutative concatenation of planar rooted trees, denoted by
$T_1\sqcup \cdots \sqcup T_b$ with $T_1,\cdots, T_b\in\calt$.
$b=\bread(F)$ is called the breadth of $F$.
The depth $\depth(F)$ of $F$ is the maximal depth of the trees $T_i, 1\leq i\leq b$.
Let $\calf$ be the set of {planar rooted forests}. Then
$\calf$ is the free semigroup generated by $\calt$ with the product $\sqcup$, and
$\bfk\,\calf$ with the product $\sqcup$ is the free noncommutative nonunitary $\bfk$-algebra on the alphabet set $\calt$.
We are going to define, for each fixed $\lambda\in\bfk$,
another product $\shpr=\shpr_\lambda$ on $\bfk\,\calf$,
making it into a unitary Rota--Baxter algebra (of weight $\lambda$).
We will suppress $\lambda$ to ease notation.

For the rest of this paper, a tree or forest means a planar rooted tree or a planar rooted forest unless otherwise specified.
Let $\lc T_1\sqcup \cdots \sqcup T_b\rc$ denote the usual {\bf grafting} of the trees $T_1,\cdots,T_b$
by adding a new root together with an edge from the new root to the root of each of the trees $T_1,\cdots, T_b$.
Let $\calf_n$, $n\geq 0$, be the set of planar rooted forests with depth less or equal to $n$. Then we have
the depth filtration
$ \calf_0 \subseteq \calf_1 \subseteq \cdots$ such that
$\calf=\cup_{n\in \NN} \calf_n$.

By using the grafting and the filtration $\calf_n$, we recursively defined in~\mcite{E-G0} a map
$$
  \shpr: \calf \times \calf \to \bfk\, \calf
$$
with the following properties
\begin{enumerate}
\item
For trees $F$ and $F'$,
\begin{equation}
F \shpr F' =\left \{ \begin{array}{ll}
    F, & {\rm\ if\ } F'=\onetree,\\
    F', & {\rm\ if\ } F=\onetree,\\
    \lc \oF\shpr \lc \oF' \rc \rc
    +\lc \lc \oF \rc \shpr \oF'\rc
    +\lambda\lc \oF \shpr \oF'\rc,
    & {\rm\ if\ } F=\lc \oF\rc, F'=\lc \oF'\rc,
    \end{array} \right .
\mlabel{eq:shprt1}
\end{equation}
\item
For forests $F=T_1\sqcup\cdots\sqcup T_b$ and
$F'=T'_1\sqcup \cdots\sqcup T'_{b'}$,
\begin{equation}
F \shpr F'= T_1\sqcup\cdots \sqcup T_{b-1}\,
\sqcup\,(T_b\shpr T'_1)\,\sqcup \,
    T'_{2}\,\cdots\,\sqcup T_{b'}.
\mlabel{eq:shprt2}
\end{equation}
\end{enumerate}
Then $\shpr$ extends to a binary operation $\shpr$ on $\bfk\calf$ by bilinearity.
As an example, we have
\begin{equation}
\td31 \shpr \tb2 = \lc \ta1\sqcup \ta1\rc \shpr \lc \ta1 \rc
= \lc (\ta1 \sqcup \ta1) \shpr \lc \ta1\rc\rc
+ \lc \lc \ta1 \sqcup \ta1\rc \shpr \ta1\rc
+ \lambda \lc (\ta1\sqcup\ta1) \shpr \ta1\rc
= \tg42 + \thj44 + \lambda \td31.
\mlabel{eq:treeex}
\end{equation}
It was shown in~\mcite{E-G0} that $(\bfk\calf,\shpr)$ is a Rota-Baxter $\bfk$-algebra.

\subsubsection{Free Rota-Baxter algebra on a set $X$}
Let $X$ be a non-empty set.
Let $F\in \calf$ with $\leaf=\leaf(F)$ leaves.
Let $X^F$ denote the set of pairs
$(F;\vec{x})$ where $\vec{x}$ is in $X^{(\ell(F)-1)}$.
Then $(F;\vec{x})$ can be identified with the forest $F$ together with an ordered decoration of $\vec{x}$ on the angles of $F$.
We use the convention that $X^{\onetree}=\{(\onetree;1)\}$.
For example, we have
$$
    \big( {\scalebox{1.15}{\tg42}}\ ;\ x     \big) = \begin{array}{l}\\[-.7cm] \xldec41r \end{array}, \quad
    \big( {\scalebox{1.15}{\thII43}}\ ;\ (x, y) \big) = \begin{array}{l}\\[-.5cm] \xyldec43 \end{array}, \quad
    \big( \ta1\ \sqcup\ {\scalebox{1.1}{\td31}}\ ;\ (x, y) \big)=  \ta1\, {x} \begin{array}{l}\\[-.3cm]
    \yldec31 \end{array}.
$$
$\ta1\, {x} \begin{array}{l}\\[-.3cm]
    \yldec31 \end{array}$ is denoted by $\ta1 \sqcup_{x} \begin{array}{l}\\[-.3cm]
    \yldec31 \end{array}$ in~\cite{E-G0}.

Let $(F;\vec{x})\in X^F$.
Let $F=T_1\sqcup \cdots \sqcup T_b$ be the decomposition of $F$
into trees. We consider the corresponding decomposition of
decorated forests. If $b=1$, then $F$ is a tree and $(F;\vec{x})$
has no further decompositions. If $b>1$,
denote $\leaf_i=\leaf(T_i), 1\leq i\leq b$.
Then
$$(T_1;(x_1,\cdots, x_{\leaf_1-1})),\
(T_2; (x_{\leaf_1+1}, \cdots, x_{\leaf_1+\leaf_2-1})),
\cdots,
(T_b; (x_{\leaf_1+\cdots+\leaf_{b-1}+1}, \cdots, x_{\leaf_1+\cdots+\leaf_b}))
$$
are well-defined angularly decorated trees when  $\leaf(T_i)>1$.
If $\leaf(T_i)=1$, then $x_{\leaf_{i-1}+\leaf_i-1}=x_{\leaf_{i-1}}$
and we use
the convention $(T_i;x_{\leaf_{i-1}+\leaf_i-1})=(T_i;\bfone)$.
With this convention, we have,
\begin{eqnarray*}
(F;(x_1,\cdots, x_{\leaf-1}))&=&
(T_1;(x_1, \cdots, x_{\leaf_1-1})){x_{\leaf_1}}
(T_2; (x_{\leaf_1+1},\cdots, x_{\leaf_1+\leaf_2-1}))
    {x_{\leaf_1+\leaf_2}}
\\
&&\cdots {x_{\leaf_1+\cdots+\leaf_{b-1}}}
(T_b; (x_{\leaf_1+\cdots+\leaf_{b-1}+1}, \cdots, x_{\leaf_1+\cdots+\leaf_b})).
\end{eqnarray*}
We call this the {\bf standard decomposition} of $(F;\vec{x})$ and
abbreviate it as
\begin{eqnarray}
(F;\vec{x})&=&(T_1;\vec{x}_1){x_{i_1}}
(T_2;\vec{x}_2) {x_{i_2}} \cdots {x_{i_{b-1}}}
(T_b;\vec{x}_b)
= D_1 {x_{i_1}} D_2 {x_{i_2}} \cdots {x_{i_{b-1}}} D_b
\mlabel{eq:stdecm}
\end{eqnarray}
where $D_i=(T_i;\vec{x}_i), 1\leq i\leq b$.
For example,
$$
    \big(\ta1\sqcup {\scalebox{1.15}{\tg42}} \sqcup {\scalebox{1.15}{\td31}}; (v, x, w, y) \big)
    = \big(\ta1;\bfone \big)\, v \big({\scalebox{1.15}{\tg42}};x)\, w \big({\scalebox{1.15}{\td31}};y \big)
    = \ta1 \,{v} \begin{array}{l}\\[-.7cm] \xldec41r \end{array} \,{w} \begin{array}{l}\\[-.3cm] \yldec31 \end{array}
$$

Let $\bfk^{NC}[X]=\bigoplus_{n\geq 0} \bfk\,X^{n}$ be the noncommutative polynomial algebra on $X$.
Denote its basis elements by vectors and its product by vector concatenation: for $\vec{x}=(x_1,\cdots,x_m),
\vec{x}'=(x'_1,\cdots,x'_n)$, define
$$(\vec{x},\vec{x}')=(x_1,\cdots,x_m,x'_1,\cdots,x'_n).$$

Define the $\bfk$-module
$$ \ncsha(X)= \bigoplus_{F\in\, \calf} \bfk\, X^{F}.$$
For $D=(F;\vec{x})\in X^F$ and
$D'=(F';\vec{x}')\in X^{F'}$, define
\begin{equation}
D\shprm D' = (F\shpr F'; (\vec{x}, \vec{x}')),
\mlabel{eq:shprm1}
\end{equation}
where $\shpr$ is defined in Eq.~(\mref{eq:shprt1}) and
Eq.~(\mref{eq:shprt2}).
For example, from Eq.~(\mref{eq:treeex}) we have
\begin{equation}
\begin{array}{l}\\[-.3cm] \xtd31 \end{array} \, \shprm \ \tb2
= \begin{array}{l}\\[-.7cm] \xldec41r \end{array}
+ \begin{array}{l}\\[-.7cm] \xthj44 \end{array}
+ \lambda \begin{array}{l}\\[-.3cm] \xtd31 \end{array}.
\mlabel{eq:treexdec}
\end{equation}

Extending the product $\shprm$ biadditively, we obtain
a binary operation
$$
\shprm:  \ncsha (X)\otimes \ncsha(X) \to \ncsha(X).
$$
For $(F;\vec{x}) \in X^F$, define
\begin{equation}
 P_X(F; \vec{x})=\lc (F;\vec{x})\rc=(\lc F \rc\, ; \vec{x})\in X^{\lc F\rc},
 \mlabel{eq:RBopm}
\end{equation}
extending to a linear operator on $\ncsha(X)$. Let
\begin{equation}
  j_X: X \to \ncsha(X)
  \mlabel{eq:jm}
\end{equation}
be the map sending $a\in X$ to $(\onetree \sqcup \onetree;a)$. The following theorem is proved in~\mcite{E-G0}.

\begin{theorem}
    The quadruple $(\ncsha(X),\shprm,P_X,j_X)$ is the free unitary Rota--Baxter algebra of
    weight $\lambda$ on the set $X$. More precisely, for any unitary Rota--Baxter algebra $(R,P)$ and map $f:X\to R$, there is a unique unitary Rota--Baxter algebra morphism
    $\free{f}: \ncsha(X) \to R$ such that $f=\free{f}\circ j_X.$
\mlabel{thm:freem}
\end{theorem}

\subsection{Free noncommutative differential Rota-Baxter algebras}
The following is the noncommutative analog of Theorem~\mref{thm:commdiffrb}.
\begin{theorem}
Let $(\bfk^{NC}\diffa{X},d^{NC}_X)
=(\bfk^{NC}[\diffs{X}],d^{NC}_X)$ be the free noncommutative differential algebra of weight $\lambda$ on a set $X$, constructed in Theorem~\mref{thm:freediff}.
Let $\ncsha(\diffs{X})$ be the free noncommutative Rota-Baxter algebra of weight $\lambda$
on $\diffs{X}$, constructed in Theorem~\mref{thm:freem}.
\begin{enumerate}
\item
There is a unique extension $\free{d}^{NC}_X$ of $d^{NC}_X$ to
$\ncsha(\diffs{X})$ so that $(\ncsha(\diffs{X}),\free{d}^{NC}_X,P_{\diffs{X}})$ is a differential Rota-Baxter algebra of weight $\lambda$.
\mlabel{it:ncdrb}
\item
The differential Rota-Baxter algebra $\ncsha(\diffs{X})$ thus obtained is the free differential Rota-Baxter algebra of weight $\lambda$ over $X$.
\mlabel{it:freencdrb}
\end{enumerate}
\mlabel{thm:ncdiffrb}
\end{theorem}
\begin{proof}
(\mref{it:ncdrb}).
We define a $\lambda$-derivation $\free{d}^{NC}_X$ on $\ncsha(\diffs{X})$ as follows. Let $F\in \calf$ and let $D\in (\diffs{X})^F$ be the forest $F$ with angular decoration by $\vec{y}\in (\diffs{X})^{\leaf(F)-1}$. Let
$$
D=(F;\vec{y})
=(T_1;\vec{y}_1){y_{i_1}}
(T_2;\vec{y}_2) {y_{i_2}} \cdots {y_{i_{b-1}}}
(T_b;\vec{y}_b)
$$
be the standard decomposition of $D$ in Eq.~(\mref{eq:stdecm}).
We define $\free{d}^{NC}_X$ by induction on the breadth $b=b(F)$ of $F$. If $b=1$, then $F$ is a tree so either $F=\onetree$ or $F=\lc \oF \rc$ for a forest $\oF$. Accordingly we define
\begin{equation}
\free{d}^{NC}_X(F;\vec{y})=\left \{\begin{array}{ll}
    0, & {\rm if\ } F=\onetree, \\
    (\oF;\vec{y}), & {\rm if\ } F=\lc \oF\rc
    \end{array} \right .
\mlabel{eq:difftree}
\end{equation}
We note that this is the only way to define $\free{d}^{NC}_X$ in order to obtain a differential Rota-Baxter algebra since $\onetree$ is the identity and $(F;\vec{y})=\lc (\oF;\vec{y})\rc$.

If $b>1$, then $F=T_1\sqcup F_t$ for another forest
$F_t=T_2\sqcup \cdots \sqcup F_b$ (t in $F_t$ stands for the tail).
So
$$
D=(F;\vec{y})=(T_1;\vec{y}_1)y_{i_1}(F_t;\vec{y}_t)
=D_1y_{i_1}D_t$$
where $D_1=(T_1;\vec{y}_1)$ and
$D_t=(T_2;\vec{y}_2) {y_{i_2}} \cdots {y_{i_{b-1}}}
(T_b;\vec{y}_b).$
We then define
\begin{eqnarray}
\free{d}^{NC}_X(D)&=&\free{d}^{NC}_X(T_1;\vec{y}_1) y_{i_1} (F_t;\vec{y}_t)
+(T_1;\vec{y}_1) d(y_{i_1}) (F_t;\vec{y}_t)
+(T_1;\vec{y}_1) y_{i_1} \free{d}^{NC}_X(F_t;\vec{y}_t) \notag\\
&&+\lambda \big(\free{d}^{NC}_X(T_1;\vec{y}_1) d(y_{i_1}) (F_t;\vec{y}_t)
+ \free{d}^{NC}_X(T_1;\vec{y}_1) y_{i_1} \free{d}^{NC}_X(F_t;\vec{y}_t)
\mlabel{eq:diffind}\\
&&+ (T_1;\vec{y}_1) d(y_{i_1}) \free{d}^{NC}_X(F_t;\vec{y}_t)\big) + \lambda^2 \free{d}^{NC}_X(T_1;\vec{y}_1) d(y_{i_1}) \free{d}^{NC}_X(F_t;\vec{y}_t),
\notag
\end{eqnarray}
where $\free{d}^{NC}_X(T_1;\vec{y}_1)$ is defined in Eq.~(\mref{eq:difftree}) and $\free{d}^{NC}_X(F_t;\vec{y}_t)$ is defined by the induction hypothesis.
Note that by Eq.~(\mref{eq:shprt2}),
$$
(T_1;\vec{y}_1)y_{i_1}(F_t;\vec{y}_t)=
(T_1;\vec{y}_1)\shpr (\onetree y_{i_1}\onetree) \shpr (F_t;\vec{y}_t).
$$
So if $\free{d}^{NC}_X$ were to satisfy the $\lambda$-Leibniz rule Eq.~(\mref{eq:diff}) with respect to the product $\shprm$, then we must have
\begin{eqnarray}
\free{d}^{NC}_X(D)&=&\free{d}^{NC}_X(T_1;\vec{y}_1)\shpr (\onetree y_{i_1} \onetree) \shpr (F_t;\vec{y}_t)
+(T_1;\vec{y}_1)\shpr \free{d}^{NC}_X(\onetree y_{i_1} \onetree) \shpr (F_t;\vec{y}_t) \notag \\
&& +\, (T_1;\vec{y}_1)  \shpr (\onetree y_{i_1} \onetree)  \shpr \free{d}^{NC}_X(F_t;\vec{y}_t)
+\lambda \free{d}^{NC}_X(T_1;\vec{y}_1) \shpr (\onetree  d^{NC}_X(y_{i_1}) \onetree) \shpr (F_t;\vec{y}_t) \notag \\
&& +\, \lambda \free{d}^{NC}_X(T_1;\vec{y}_1) \shpr (\onetree  y_{i_1} \onetree)  \shpr \free{d}^{NC}_X(F_t;\vec{y}_t)
+ \lambda (T_1;\vec{y}_1)  \shpr (\onetree d^{NC}_X(y_{i_1})  \onetree)  \shpr \free{d}^{NC}_X(F_t;\vec{y}_t)
\mlabel{eq:diffind2}\\
&& +\, \lambda^2 \free{d}^{NC}_X(T_1;\vec{y}_1)  \shpr (\onetree d^{NC}_X(y_{i_1})  \onetree) \shpr \free{d}^{NC}_X(F_t;\vec{y}_t).
\notag
\end{eqnarray}
Since $\free{d}^{NC}_X$ is to extend $d_X^{NC}:\bfk^{NC}\diffa{X} \to
\bfk^{NC} \diffa{X}$, we have
$$ \free{d}^{NC}_X(\onetree\, y_{i_1} \onetree)
=\free{d}^{NC}_X(j_{\diffs{X}}(y_{i_1}))= j_{\diffs{X}}(d^{NC}_X(y_{i_1}))=\onetree \, d^{NC}_X(y_{i_1}) \onetree.$$
So by Eq.~(\mref{eq:shprt2}), Eq.~(\mref{eq:diffind2}) agrees with Eq.~(\mref{eq:diffind}). Thus
$\free{d}^{NC}_X(D)$ is the unique map that satisfies the $\lambda$-Leibniz rule~(\mref{eq:diff}).

We also have the short hand notation,
\begin{equation}
\free{d}^{NC}_X(D)=\free{d}^{NC}_X(D_1)y_{i_1}D_t+D_1 \free{d}^{NC}_X(y_{i_1}D_t) +\lambda \free{d}^{NC}_X(D_1)\free{d}^{NC}_X(y_{i_1}D_t),
\mlabel{eq:difft}
\end{equation}
where
$$\free{d}^{NC}_X(y_{i_1}D_t):= d^{NC}_X(y_{i_1})D_t +y_{i_1}\free{d}^{NC}_X(D_t) +\lambda d^{NC}_X(y_{i_1})\free{d}^{NC}_X(D_t).$$
Similarly, we can also write $D=D_h y_{i_{b-1}} D_{b}$
where $D_h$ ($h$ stands for the head) is a angularly decorated forest and $D_{b}$ is a angularly decorated tree. Then
\begin{equation}
\free{d}^{NC}_X(D)=\free{d}^{NC}_X(D_hy_{i_{b-1}})D_{b}+D_hy_{i_{b-1}} \free{d}^{NC}_X(D_{b})+\lambda \free{d}^{NC}_X(D_hy_{i_{b-1}})\free{d}^{NC}_X(D_{b}).
\mlabel{eq:diffh}
\end{equation}
In fact, write
$$D=v_1v_2\cdots v_{2b-1},$$
where
$$v_j=\left \{\begin{array}{ll}
D_{(j-1)/2}, & j {\rm\ odd}, \\
y_{i_{j/2}}, & j {\rm\ even}.
\end{array} \right . $$
Then using Eq.~(\mref{eq:diffind}) and an induction on $b$, we obtain the ``general Leibniz formula" of weight $\lambda$ with respect to the concatenation product:
\begin{equation}
\free{d}^{NC}_X(D)=\sum_{I\subseteq [2b-1]} \lambda^{|I|-1} v_{I,1}v_{I,2}\cdots v_{I,2b-1},
\mlabel{eq:diffgen}
\end{equation}
where $[2b-1]=\{1,\cdots,2b-1\}$ and
$$v_{I,j} =\left\{ \begin{array}{ll}
    v_j,  j\not\in I, \vspace{.2cm}\\
    \free{d}^{NC}_X(v_j), j\in I, j {\rm\ odd},
    \vspace{.2cm}\\
    d^{NC}_X(v_j), j\in I, j {\rm\ even}.
    \end{array} \right .
$$

We now prove that $\free{d}^{NC}_X$ is a derivation of weight $\lambda$ with respect to the product $\shprm$. Let $D$ and $D'$ be angularly decorated forests and write
$$D=(F;\vec{y})=(T_1;\vec{y}_1){y_{i_1}}
(T_2;\vec{y}_2) {y_{i_2}} \cdots {y_{i_{b-1}}}
(T_{b};\vec{y}_{b})
=D_h y_{i_{b-1}} D_b$$
 and
$$D'
 =(F';\vec{y}')=(T'_1;\vec{y}'_1){y'_{i_1}}
(T'_2;\vec{y}'_2) {y'_{i_2}} \cdots {y'_{i_{b'-1}}}
(T'_{b'};\vec{y}'_{b'})
 =D'_1 y'_{i_1}D'_t$$
be as above with angularly decorated trees $D_b$, $D'_1$, angularly decorated forests $D_h$, $D'_t$ and $y_{i_{b-1}},y'_{i_1}\in \diffs{X}$.
Then by Eq.~(\mref{eq:shprt2}) (see~\mcite{E-G0} for further details), $D \shprm D'$ has the standard decomposition
\begin{eqnarray}
D \shprm D' &=&(T_1;\vec{y}_1){y_{i_1}} \cdots {y_{i_{b-1}}}
\big((T_b;\vec{y}_b)\shprm (T'_1;\vec{y}'_1)\big) {y'_{i_1}}
 \cdots {y'_{i_{b'-1}}} (T'_{b'};\vec{y}'_{b'}) \notag\\
 &=& D_h y_{i_{b-1}} (D_b\shpr D'_1) y'_{i_1}D'_t
\mlabel{eq:shprm2}
\end{eqnarray}
where
\begin{eqnarray}
\lefteqn{D_b\shpr D'_1=(T_b;\vec{y}_b)\shprm (T'_1;\vec{y}'_1)} \mlabel{eq:shprm3}\\
&=& \left \{ \begin{array}{ll}
(\onetree; \bfone), & {\rm if\ } T_b=T'_1=\onetree \ ({\rm so\ }
                \vec{y}_b=\vec{y}'_1=\bfone),\\
(T_b,\vec{y}_b), & {\rm if\ } T'_1=\onetree, T_b\neq \onetree, \\
(T'_1,\vec{y}'_1), & {\rm if\ } T'_1\neq\onetree, T_b= \onetree, \\
\lc (T_b;\vec{y})\shprm (\oF'_1;\vec{y}')\rc + \lc (\oF_b;\vec{y})\shprm (T'_1;\vec{y}') \rc  & \\
+\lambda \lc (\oF_b;\vec{y})\shprm (\oF'_1;\vec{y}')\rc, &
    {\rm if\ } T'_1=\lc \oF'_1\rc \neq\onetree, T_b=\lc \oF_b\rc \neq \onetree.
\end{array} \right .
\notag
\end{eqnarray}

By Eq.~(\mref{eq:shprm2}) and Eq.~(\mref{eq:diffgen}), we have
\begin{eqnarray}
\free{d}^{NC}_X(D\shprm D')&=&
\free{d}^{NC}_X\big((D_h y_{i_{b-1}})(D_b\shprm D'_1)(y'_{i_1}D'_{b'})\big) \notag\\
&=& \free{d}^{NC}_X(D_h y_{i_{b-1}})(D_b\shprm D'_1)(y'_{i_1}D'_{b'})
+ (D_h y_{i_{b-1}})\free{d}^{NC}_X(D_b\shprm D'_1)(y'_{i_1}D'_{b'})
\notag \\
&&+ (D_h y_{i_{b-1}})(D_b\shprm D'_1)\free{d}^{NC}_X(y'_{i_1}D'_{b'})
+\lambda \free{d}^{NC}_X(D_h y_{i_{b-1}})\free{d}^{NC}_X(D_b\shprm D'_1) (y'_{i_1}D'_{b'}) \mlabel{eq:diffprod}\\
&& + \lambda\, \free{d}^{NC}_X(D_h y_{i_{b-1}})(D_b\shprm D'_1) \free{d}^{NC}_X(y'_{i_1}D'_{b'})
+ \lambda\, (D_h y_{i_{b-1}})\free{d}^{NC}_X(D_b\shprm D'_1) \free{d}^{NC}_X(y'_{i_1}D'_{b'}) \notag\\
&&+\lambda^2
\free{d}^{NC}_X(D_h y_{i_{b-1}})\free{d}^{NC}_X(D_b\shprm D'_1) \free{d}^{NC}_X(y'_{i_1}D'_{b'}). \notag
\end{eqnarray}
Using Eq.~(\mref{eq:shprm3}), we have
\begin{equation}
\free{d}^{NC}_X(D_b\shprm D'_1)=
\free{d}^{NC}_X (D_b)\shprm D'_1 +D_b\shprm \free{d}^{NC}_X(D'_1)
+\lambda \free{d}^{NC}_X(D_b)\shprm \free{d}^{NC}_X(D'_1).
\mlabel{eq:diffirr}
\end{equation}
Applying this to Eq.~(\mref{eq:diffprod}), we find that
the resulting expansion for $\free{d}^{NC}_X(D\shprm D')$ agrees with the expansion of
$$ \free{d}^{NC}_X(D) \shprm D'+ D\shprm \free{d}^{NC}_X(D')
+ \lambda \free{d}^{NC}_X(D) \shprm \free{d}^{NC}_X(D')$$
after applying Eq.~(\mref{eq:difft}) to $\free{d}^{NC}_X(D)$ and applying Eq.~(\mref{eq:diffh}) to $\free{d}^{NC}_X(D')$.

As an example, from Eq.~(\mref{eq:treexdec}), we have
\begin{equation}
\free{d}^{NC}_X(\begin{array}{l}\\[-.3cm] \xtd31 \end{array} \, \shprm \ \tb2)
= \free{d}^{NC}_X\big( \begin{array}{l}\\[-.7cm] \xldec41r \end{array}
+ \begin{array}{l}\\[-.7cm] \xthj44 \end{array}
+ \lambda \begin{array}{l}\\[-.3cm] \xtd31 \end{array}\big)
= \onetree\, x \tb2
+ \begin{array}{l}\\[-.3cm] \xtd31 \end{array}
+ \onetree\, x \onetree\,.
\mlabel{eq:treexdecd}
\end{equation}
This agrees with
$$
\free{d}^{NC}_X(\begin{array}{l}\\[-.3cm] \xtd31 \end{array}) \, \shprm \ \tb2
+ \begin{array}{l}\\[-.3cm] \xtd31 \end{array} \, \shprm \ \free{d}^{NC}_X(\tb2)
+ \lambda \free{d}^{NC}_X(\begin{array}{l}\\[-.3cm] \xtd31 \end{array}) \, \shprm \ \free{d}^{NC}_X(\tb2).$$

\smallskip

\noindent
(\mref{it:freencdrb}).
The proof of the freeness of $\ncsha(D(X))$ as a free differential Rota-Baxter algebra of weight $\lambda$ is the same as the proof of the freeness of $\sha(D(X))$ in Theorem~\mref{thm:commdiffrb}.
\end{proof}

\section{Structure of a differential algebra on forests}
\mlabel{sec:tree}
We now give the structure of a differential Rota-Baxter algebra of weight $\lambda$ to rooted forests without decorations. It should be possible to derive this as a special case from a suitable generalization of the construction in Theorem~\mref{thm:ncdiffrb}. To avoid making the process too complicated, we give a direct construction. See~\mcite{G-L} for the work of Grossman and Larson on differential algebra structures on their Hope algebra of trees.

Let $(\bfk \calf, \shpr, \lc\ \rc)$ be the Rota-Baxter algebra of planar rooted forests defined in Section~\mref{sss:rbtree}.
Let $F\in \calf$ be a rooted forest. By Eq.~(\mref{eq:shprt2}), the unique decomposition $F=T_1\sqcup \cdots \sqcup T_b$ into
rooted trees $T_1,\cdots,T_b\in \calt$ gives
the decomposition
\begin{equation}
F= T_1 \shpr (\onetree\sqcup \onetree) \shpr \cdots
\shpr (\onetree \sqcup \onetree) \shpr T_b.
\mlabel{eq:treedecomp1}
\end{equation}
Denote this by
\begin{equation}
F=V_1 \shpr V_2 \shpr \cdots \shpr V_{2b-1},
\mlabel{eq:treedecomp2}
\end{equation}
where
$$ V_i=\left \{\begin{array}{ll}
    T_{(i+1)/2}, & i {\rm\ odd},\\
    (\onetree\sqcup\onetree), & i {\rm\ even}
    \end{array} \right .
$$
We call Eq.~(\mref{eq:treedecomp2}) the {\bf $\shpr$-standard decomposition} of $F$.
This decomposition is unique since it is uniquely determined by the unique decomposition of $F$ into rooted trees.

We define a linear operator
\begin{equation}
d_\calf: \bfk \calf \to \bfk \calf.
\mlabel{eq:treediff1}
\end{equation}
as follows. First let $V$ be either $\onetree \sqcup \onetree$ or a tree, hence of the form $\onetree$ or $\lc \oV \rc$ for a forest $\oV$. Define
\begin{equation}
d_\calf(V) = \left \{\begin{array}{ll}
0, & V= \onetree, \\
1, & V=\onetree \sqcup \onetree, \\
\oV, & V=\lc \oV \rc.
\end{array} \right .
\mlabel{eq:treediff2}
\end{equation}
Next let $F\in \calf$ have the $\shpr$-standard decomposition in Eq.~(\mref{eq:treedecomp2}). Define
\begin{equation}
d_\calf(F)=\sum_{\emptyset \neq I\subseteq [k]}
    \lambda^{|I|-1} V_{I,1} \shpr \cdots \shpr V_{I,k},
\mlabel{eq:treediff3}
\end{equation}
where for $I\subseteq [k]$,
\begin{equation}
V_{I,i}=\left \{ \begin{array}{ll}
    V_i, & i\not\in I, \\
    d_\calf(V_i), & i\in I
\end{array} \right .
\mlabel{eq:treediff4}
\end{equation}
with $d_\calf(V_i)$ as defined in Eq.~(\mref{eq:treediff2}).
Finally extend $d_\calf$ to $\bfk \calf$ by $\bfk$-linearity.

It is clear that $d_\calf$ satisfies the recursive relation
\begin{equation}
d_\calf(F) = d_\calf(V_1)\shpr (V_2\shpr \cdots \shpr V_k)
+ V_1 \shpr d_\calf(V_2\shpr \cdots \shpr V_k)
+ \lambda d_\calf(V_1) \shpr d_\calf(V_2\shpr \cdots \shpr V_k).
\mlabel{eq:treediff5}
\end{equation}

We give some examples. By the third case in Eq.~(\mref{eq:treediff2}), we have
\begin{equation}
d_\calf(\td31)=\onetree\sqcup \onetree,\quad
d_\calf(\tg42)= \onetree\sqcup \tb2.
\end{equation}
Further, since $\onetree \sqcup \tb2=
   (\onetree\sqcup \onetree) \shpr \tb2$, we have
\begin{equation}
d_\calf(\onetree \sqcup \tb2)
    = \onetree\shpr \tb2 + (\onetree\sqcup\onetree) \shpr \onetree +\lambda (\onetree \shpr \onetree)\\
    = \tb2 + \onetree\sqcup \onetree + \lambda\, \onetree.
    \mlabel{eq:treeex1}
\end{equation}
Similarly,
\begin{eqnarray}
d_\calf(\td31 \sqcup \onetree)&=&
d_\calf(\td31 \shpr (\onetree \sqcup \onetree))
\notag\\
&=& (\onetree\sqcup \onetree)\shpr (\onetree\sqcup\onetree) + \td31 \shpr \onetree + \lambda (\onetree\sqcup \onetree)\shpr \onetree
\mlabel{eq:treeex2}\\
&=& \onetree\sqcup\onetree\sqcup\onetree
    + \td31 +  \lambda (\onetree\sqcup \onetree).
\notag
\end{eqnarray}
As another example, from the $\shpr$-standard decomposition
$$\td31\sqcup \tb2=
    \td31 \shpr (\onetree\sqcup\onetree) \shpr \tb2,
$$
by Eq.~(\mref{eq:treeex}), Eq.~(\mref{eq:treediff5}) and Eq.~(\mref{eq:treeex1}) we have
\begin{eqnarray*}
d_\calf(\td31\sqcup \tb2)&=&
d_\calf(\td31)\shpr \big((\onetree\sqcup\onetree) \shpr \tb2\big)
    + \td31 \shpr d_\calf\big((\onetree\sqcup\onetree) \shpr \tb2\big)
    + \lambda d_\calf(\td31) \shpr d_\calf\big((\onetree\sqcup\onetree) \shpr \tb2\big)\\
&=& (\onetree\sqcup\onetree)\shpr \big((\onetree\sqcup\onetree) \shpr \tb2\big)
+ \td31 \shpr \big(\tb2 + \onetree\sqcup \onetree + \lambda\, \onetree\big)
+ \lambda (\onetree\sqcup \onetree) \shpr (\tb2 + \onetree\sqcup \onetree + \lambda\, \onetree)\\
&=&
\onetree\sqcup\onetree\sqcup \onetree \sqcup \tb2 + \td31\shpr \tb2 +\td31\sqcup \onetree + \lambda \td31+\lambda \onetree\sqcup \tb2 +\lambda \onetree\sqcup\onetree\sqcup\onetree +\lambda^2 \onetree\sqcup\onetree\\
&=& \tg42 +\thj44+2\lambda \td31 +\td31\sqcup \onetree
+\onetree\sqcup\onetree\sqcup\onetree\sqcup\tb2 +\lambda \onetree\sqcup\tb2
+\lambda \onetree\sqcup\onetree\sqcup\onetree
+\lambda^2 \onetree\sqcup\onetree.
\end{eqnarray*}

\begin{theorem}
The triple $(\bfk \calf, d_\calf, \lc\ \rc)$ is a differential Rota-Baxter algebra of weight $\lambda$.
\mlabel{thm:treediff}
\end{theorem}
\begin{proof}
By the third case of Eq.~(\mref{eq:treediff2}), $d_\calf \circ \lc\ \rc=\id$. So we only need to show that $d_\calf$ is a differential operator of weight $\lambda$, that is, $d_\calf$ satisfies the $\lambda$-Leibniz rule in Eq.~(\mref{eq:diff}):
\begin{equation}
d_\calf(F\shpr F')=d_\calf(F)\shpr F' + F\shpr d_\calf(F') + \lambda d_\calf(F) \shpr d_\calf(F').
\mlabel{eq:diff2}
\end{equation}
This is not immediate since the $\shpr$-standard decomposition of $F\shpr F'$ is not the product of
the $\shpr$-standard decomposition of $F$ and $F'$.

First let $F$ and $F'$ be trees. Then $F$ is either $\onetree$ or $\lc \oF\rc$ for a forest $\oF$. Similarly for $F'$. Since $\onetree$ is the unit, Eq.~(\mref{eq:diff}) trivially holds if $F=\onetree$ or $F'=\onetree$. If $F=\lc \oF\rc$ and $F'=\lc \oF' \rc$.
Then by the Rota-Baxter equation (\mref{eq:rba}) and Eq.~(\mref{eq:treediff4}), we have
\begin{equation}
 d_\calf(F\shpr F')= \oF \shpr F' + F\shpr \oF' +\lambda \oF \shpr \oF'.
 \mlabel{eq:diffprod0}
\end{equation}
This is Eq.~(\mref{eq:diff2}).

In general, let $F$ and $F'$ be forests and let
$$F=V_1\shpr \cdots \shpr V_{2b-1},
\qquad F'=V'_1\shpr \cdots \shpr V'_{2b'-1}$$
be their $\shpr$-standard decompositions from
Eq.~(\mref{eq:treedecomp1}). Then
$$ F\shpr F'= V_1\shpr \cdots \shpr V_{2b-2} \shpr (V_{2b-1} \shpr V'_1) \shpr V'_2\shpr \cdots \shpr V'_{2b'-1}
$$
is the $\shpr$-standard decomposition of $F\shpr F'$.
Here $V_{2b-1}\shpr V'_1=\sum_k Z\,''_k$ is a tree or a linear combination of trees $Z\,''_k$ given in Eq.~(\mref{eq:shprt1}).
As in Eq.~(\mref{eq:treedecomp2}), we rewrite it as
$$F\shpr F'= W_1\shpr \cdots \shpr W_{2(b+b'-1)-1}.$$
In particular,
$ W_{2b-1}=V_{2b-1}\shpr V'_1=\sum_k Z\,''_k$.
Then by definition,
\begin{equation}
 d_\calf(F\shpr F') = \sum_{\emptyset \neq J\subseteq [2(b+b'-1)-1]} \lambda^{|J|-1} W_{J,1}\shpr \cdots \shpr W_{J,2(b+b'-1)-1}
\mlabel{eq:diffprod1}
\end{equation}
with
$W_{J,j}$ defined in the same way as $T_{I,i}$ in Eq.~(\mref{eq:treediff2}) and
$d_\calf(W_{2b-1})=\sum_k d_\calf(Z\,''_k).$
Depending on whether or not $2b-1\in J$, we can rewrite
Eq.~(\mref{eq:diffprod1}) as
{\small
\begin{eqnarray}
d_\calf(F\shpr F') &=&
\sum_{2b-1\in J\subseteq [2(b+b'-1)-1]}
\lambda^{|J|-1} W_{J,1}\shpr \cdots \shpr W_{J,2b-2}\shpr d_\calf(W_{2b-1}) \shpr \cdots \shpr W_{J,2(b+b'-1)-1} \notag\\
&& + \sum_{2b-1\not\in J\subseteq [2(b+b'-1)-1]}
\lambda^{|J|-1} W_{J,1}\shpr \cdots \shpr W_{J,2b-2}\shpr W_{J,2b-1} \shpr \cdots \shpr W_{J,2(b+b'-1)-1} \notag \\
&=& \Big(\sum_{\check{I}\subseteq [2b-2]}
    \lambda^{|\check{I}|}(V_{\check{I},1} \shpr \cdots \shpr V_{\check{I},2b-2}\Big) \shpr d_\calf(V_{2b-1}\shpr V'_1) \mlabel{eq:leib1} \\
&& \qquad \shpr \Big(\sum_{\check{I}'\subseteq \{2,\cdots,2b'-1\}} \lambda^{|\check{I}'|} V'_{\check{I}',2} \shpr \cdots
    \shpr V'_{\check{I}',2b'-1} \Big) \notag\\
&&+ \hspace{-1cm}
\sum_{\scriptsize{\begin{array}{l}
\check{I}\subseteq \{1,\cdots,2b-2\}\\
\check{I}'\subseteq \{2,\cdots,2b'-1\}\\
\check{I}\neq \emptyset {\rm\ or\ } \check{I}'\neq\emptyset\end{array} }}
\hspace{-1cm}
    \lambda^{|\check{I}|+|\check{I}'|-1} (V_{\check{I},1} \shpr \cdots \shpr V_{\check{I},{2b-2}})
    \shpr (V_{2b-1}\shpr V'_1) \shpr (V'_{\check{I}',2} \shpr \cdots
    \shpr V'_{\check{I}',2b'-1} ) \notag
\end{eqnarray}
}
By Eq.~(\mref{eq:diffprod0}),
$$ d_{\calf} (V_{2b-1}\shpr V'_1)=d_\calf (V_{2b-1})\shpr V'_1 +V_{2b-1}\shpr d_\calf(V'_1)+\lambda d_\calf (V_{2b-1})\shpr d_\calf (V'_1).$$
Denote $I\subseteq [2b-1]$ and $I'\subseteq [2b'-1]$.
We can write the first sum in Eq.~(\mref{eq:leib1}) as
\begin{eqnarray}
&& \Big(\hspace{-.3cm} \sum_{\scriptsize{\begin{array}{l} 2b-1\in I\\
1\not\in I'\end{array}}}
\hspace{-.3cm} + \sum_{\scriptsize{\begin{array}{l} 2b-1\not\in I\\
1\in I'\end{array}}}
\hspace{-.3cm} +\sum_{\scriptsize{\begin{array}{l} 2b-1\in I\\
1\in I'\end{array}}}\hspace{-.3cm} \Big)
\lambda^{|I|+|I'|-1} V_{I,1}\shpr \cdots \shpr (V_{I,2b-1}\shpr V'_{I',1})\shpr \cdots \shpr V'_{I',2b'-1}
\notag \\
&&
= \Big(\hspace{-.3cm}  \sum_{\scriptsize{\begin{array}{l}
2b-1\in I\\
I'= \emptyset\end{array}}}
\hspace{-.3cm} +
\sum_{\scriptsize{\begin{array}{l}
2b-1\in I\\
1\not\in I'\neq \emptyset \end{array}}}
\hspace{-.3cm} +
\sum_{\scriptsize{\begin{array}{l} I= \emptyset \\
1\in I'\end{array}}}
\hspace{-.3cm} + \sum_{\scriptsize{\begin{array}{l} 2b-1\not\in I\neq \emptyset\\
1\in I'\end{array}}}
\hspace{-.3cm} +\sum_{\scriptsize{\begin{array}{l} 2b-1\in I\\
1\in I'\end{array}}}\hspace{-.3cm} \Big) \mlabel{eq:leib2} \\
&& \qquad \qquad \lambda^{|I|+|I'|-1} V_{I,1}\shpr \cdots \shpr (V_{I,2b-1}\shpr V'_{I',1})\shpr \cdots \shpr V'_{I',2b'-1}.
\notag
\end{eqnarray}

For the second sum in Eq.~(\mref{eq:leib1}), we have
\begin{equation}
\Big(\hspace{-.3cm}  \sum_{\scriptsize{\begin{array}{l} 2b-1\not\in I\neq \emptyset\\ I'=\emptyset\end{array}}}
\hspace{-.3cm} +\hspace{-.3cm}  \sum_{\scriptsize{\begin{array}{l} I= \emptyset\\ 1\not\in I'\neq \emptyset\end{array}}}
\hspace{-.3cm} +\hspace{-.3cm}
\sum_{\scriptsize{\begin{array}{l} 2b-1\not\in I\neq \emptyset\\ 1\not\in I'\neq \emptyset\end{array}}}
\hspace{-.3cm} \Big)
\lambda^{|I|+|I'|-1} V_{I,1}\shpr \cdots \shpr (V_{I,2b-1}\shpr V'_{{I'},1})\shpr \cdots \shpr V'_{{I'},2b'-1}.
\mlabel{eq:leib3}
\end{equation}
The first sum on the right hand side of  Eq.~(\mref{eq:leib2}) adding to the first sum in Eq.~(\mref{eq:leib3}) gives
\begin{eqnarray*}
&& \sum_{I\neq \emptyset, I'=\emptyset}
\lambda^{|I|+|I'|-1} V_{I,1}\shpr \cdots \shpr (V_{I,2b-1}\shpr V'_1)\shpr \cdots \shpr V'_{2b'-1}\\
&=& \sum_{I\neq \emptyset} \lambda^{|I|-1}
V_{I,1}\shpr \cdots \shpr (V_{I,2b-1}\shpr V'_1)\shpr \cdots \shpr V'_{2b'-1}\\
&=& d_\calf(F)\shpr F'.
\end{eqnarray*}
Similarly, the third sum on the right hand side of Eq.~(\mref{eq:leib2}) adding to the second term in Eq.~(\mref{eq:leib3}) gives
$F \shpr d_\calf(F')$. The remaining terms on the right hand side of Eq.~(\mref{eq:leib2}) and Eq.~(\mref{eq:leib3}) add to $\lambda d_\calf(F)\shpr d_\calf(F')$. This proves the $\lambda$-Leibniz rule (\mref{eq:diff}).
\end{proof}


\end{document}